\documentclass[12pt]{amsart}
\usepackage{a4wide}
\usepackage[utf8]{inputenc}
\usepackage{amsmath}
\usepackage{amsfonts}
\usepackage{amssymb}
\usepackage{amsthm}
\usepackage{subfigure}
\usepackage{xcolor}
\newtheorem{theo}{Theorem}[section]
\newtheorem{prop}[theo]{Proposition}

\newtheorem{lemm}[theo]{Lemma}

\theoremstyle{definition}

\theoremstyle{remark}
\newtheorem{rema}[theo]{Remark}
\newcommand{\Op}{\operatorname{Op}}

\newcommand{\nwc}{\newcommand}
\nwc{\eps}{\epsilon}
\nwc{\ep}{\epsilon}
\nwc{\vareps}{\varepsilon}
\nwc{\Oph}{\operatorname{Op}_\hbar}
\nwc{\la}{\langle}
\nwc{\ra}{\rangle}

\nwc{\mf}{\mathbf} %Latex (as in \bf not tilted math letters)
\nwc{\blds}{\boldsymbol} %Latex
\nwc{\ml}{\mathcal} %Latex

\nwc{\defeq}{\stackrel{\rm{def}}{=}}

\nwc{\cE}{\ml{E}}
\nwc{\cN}{\ml{N}}
\nwc{\cO}{\ml{O}}
\nwc{\cP}{\ml{P}}
\nwc{\cU}{\ml{U}}
\nwc{\cV}{\ml{V}}
\nwc{\cW}{\ml{W}}
\nwc{\tU}{\widetilde{U}}
\nwc{\IN}{\mathbb{N}}
\nwc{\IR}{\mathbb{R}}
\nwc{\IZ}{\mathbb{Z}}
\nwc{\IC}{\mathbb{C}}
\nwc{\IT}{\mathbb{T}}
\nwc{\IS}{\mathbb{S}}
\nwc{\tP}{\widetilde{P}}
\nwc{\tPi}{\widetilde{\Pi}}
\nwc{\tV}{\widetilde{V}}
\nwc{\supp}{\operatorname{supp}}
\nwc{\rest}{\restriction}
\nwc{\todo}[1]{$\clubsuit$ {\tt #1}}

\begin{document}

\title[Local $L^p$ norms of Schr\"odinger eigenfunctions]{Local $L^p$ norms of Schr\"odinger eigenfunctions on $\mathbb{S}^2$}

\author[Gabriel Rivi\`ere]{Gabriel Rivi\`ere}

\address{Laboratoire de math\'ematiques Jean Leray (U.M.R. CNRS 6629), Universit\'e de Nantes, 2 rue de la Houssini\`ere, BP92208, 
44322 Nantes Cedex 3, France}
\address{Institut Universitaire de France}
\email{gabriel.riviere@univ-nantes.fr}

\begin{abstract}
On the canonical $2$-sphere and for Schr\"odinger eigenfunctions, we obtain a simple geometric criterion on the potential under which we can improve, near a given point and for every $p\neq 6$, Sogge's estimates by a power of the eigenvalue. This criterion can be formulated in terms of the critical points of the Radon transform of the potential and it is independent of the choice of eigenfunctions.\\

%[French version] Sur la $2$-sph\`ere munie de sa m\'etrique canonique, nous obtenons un crit\`ere g\'eom\'etrique simple sur le potentiel pour lequel nous am\'eliorons d'un facteur polynomial les estim\'ees de Sogge sur les fonctions propres de Schr\"odinger pr\`es d'un point fix\'e et pour tout $p\neq 6$. Ce crit\`ere peut \^etre formul\'e en termes des points critiques de la transform\'ee de Radon du potentiel et il est ind\'ependant du choix des fonctions propres.
\end{abstract}

\maketitle

\section{Introduction}

The purpose of this work is to study high frequency asymptotics of eigenfunctions to the Schr\"odinger operator on the $2$-sphere
\begin{equation}\label{e:sphere}\mathbb{S}^2:=\left\{(x_1,x_2,x_3)\in\mathbb{R}^3:x_1^2+x_2^2+x_3^2=1\right\}.\end{equation}
We endow $\IS^2$ with the Riemannian metric $g_0$ induced by the Euclidean metric on $\mathbb{R}^3$. In that geometric context and given an element $V\in\mathcal{C}^{\infty}(\mathbb{S}^2,\mathbb{R})$, there exists an orthonormal basis~\cite[Th.~14.7]{Zw12} of $L^2(\mathbb{S}^2,d\upsilon_{g_0})$ made of solutions to
\begin{equation}\label{e:schrodinger}-\Delta_{g_0}\psi_\lambda +V\psi_\lambda =\lambda^2\psi_\lambda,\quad\lambda^2\in\mathbb{R},\end{equation}
where $\Delta_{g_0}$ is the Laplace-Beltrami operator and $d\upsilon_{g_0}$ is the Riemannian volume, both induced by $g_0$. By elliptic regularity, solutions to~\eqref{e:schrodinger} are smooth~\cite[\S~14.3]{Zw12} and a classical Theorem of Sogge~\cite{So88} states that, for every $2\leq p\leq +\infty$, there exists $C_p>0$ such that, for any solution $(\psi_\lambda,\lambda)$ to~\eqref{e:schrodinger},
\begin{equation}\label{e:Lp-sogge}\|\psi_\lambda\|_{L^p(\IS^2)}\leq C_p(1+|\lambda|)^{\sigma_0(p)}\|\psi_\lambda\|_{L^2(\IS^2)},\end{equation}
where\footnote{The case $p=\infty$ is a consequence of the local Weyl law~\cite{Ho68}.} 
$$\sigma_0(p):=\max\left\{\frac{1}{4}-\frac{1}{2p},\frac{1}{2}-\frac{2}{p}\right\}.$$
The critical exponent for which both quantities in the maximum coincide is given by $p_c=6$. Recall that Sogge's result can be extended to $\ml{O}(\lambda)$-quasimodes of $-\Delta_{g_0}$ as solutions to~\eqref{e:schrodinger} are -- see also~\cite{KTZ}. In the case where $V\equiv 0$, these upper bounds are optimal using appropriate sequences of spherical harmonics~\cite{So15}. However, for generic sequences~\cite{VdK97, Ze08, BuLe13} or for families satisfying certain extra invariance properties~\cite{BrLM20}, these bounds can drastically be improved when $V\equiv 0$. 

Our aim is to show that the presence of a potential allows to improve~\eqref{e:Lp-sogge} away from certain critical geodesics and for \emph{any} sequence of eigenfunctions. In order to state our results, we introduce the space of oriented closed geodesics $G(\mathbb{S}^2)$ of the sphere. By identifying each oriented closed geodesic with an oriented plane of $\mathbb{R}^3$, $G(\mathbb{S}^2)$ is diffeomorphic to $\mathbb{S}^2$. Through this identification, $G(\mathbb{S}^2)\simeq\mathbb{S}^2$ is endowed with the symplectic structure induced by the one on the cotangent bundle $T^*\mathbb{S}^2$~\cite[p.~58]{Bes78}. We also define the Radon transform of the potential $V$:
$$\mathcal{R}(V):\gamma\in G(\mathbb{S}^2)\mapsto \frac{1}{2\pi}\int_0^{2\pi} V(\gamma(s))ds \in\mathbb{R},$$
which belongs to $\mathcal{C}^{\infty}(G(\mathbb{S}^2))$. Thanks to the symplectic structure on $G(\IS^2)$, one can define its Hamiltonian vector field $X_{\langle V\rangle}$. We denote its critical points by
$$\text{Crit}(\mathcal{R}(V)):=\left\{\gamma\in G(\mathbb{S}^2): D_\gamma\mathcal{R}(V)=0\right\}=\left\{\gamma\in G(\mathbb{S}^2): X_{\langle V\rangle}(\gamma)=0\right\}.$$
Observe that $\ml{R}(V)$ is always an even function on $G(\IS^2)$. In particular, it can be identified with a function on $\mathbb{R}P^2$ and it has thus at least $6$ critical points on $G(\IS^2)$ by Morse inequalities. In fact, Guillemin showed~\cite{Gu76} that
$$\mathcal{R}:\mathcal{C}^{\infty}(\mathbb{S}^2)\rightarrow \mathcal{C}^{\infty}_{\text{even}}(G(\mathbb{S}^2))\simeq\mathcal{C}^{\infty}(\IR P^2)$$
is a surjective map whose kernel is equal to $\mathcal{C}^{\infty}_{\text{odd}}(\mathbb{S}^2)$. As a corollary, for a generic choice of $V$ in the $\mathcal{C}^{\infty}$-topology, $\text{Crit}(\mathcal{R}(V))$ is a finite set. Finally, given $x_0\in\IS^2$, we set 
$$\Gamma_{x_0}:=\left\{\gamma\in G(\IS^2):\ x_0\in\gamma\right\}.$$
Our main result reads as follows
\begin{theo}\label{t:maintheo} Let $x_0\in \IS^2$ such that
\begin{equation}\label{e:hyp-crit}\operatorname{Crit}(\mathcal{R}(V))\cap \Gamma_{x_0}=\emptyset,\end{equation}
and
\begin{equation}\label{e:hyp-trans}
%\left|\left\{\gamma \in\Gamma_{x_0}:\ X_{\langle V\rangle}(\gamma)\in T_\gamma\Gamma_{x_0}\right\}\right|<\infty.
\ml{R}(V)|_{\Gamma_{x_0}}\ \text{is a Morse function.}\end{equation}
Then, there exists $r_0>0$ such that, for every $2\leq p\leq +\infty$, one can find $C_{x_0,p}>0$ so that, for any solution $(\psi_\lambda,\lambda)$ to~\eqref{e:schrodinger},
$$\|\psi_\lambda\|_{L^p(B_{r_0}(x_0))}\leq C_{x_0,p}(\log(2+|\lambda|))^{\varepsilon(p)}(1+|\lambda|)^{\sigma_0(p)-\delta(p)}\|\psi_\lambda\|_{L^2(\mathbb{S}^2)},$$
where $B_{r_0}(x_0)$ is the closed (geodesic) ball of radius $r_0$ centered at $x_0$ and where, for $4< p\leq \infty,$
$$\delta(p):=\frac{1}{18}\left|1-\frac{6}{p}\right|,\quad \varepsilon(p)=0$$
and, for $2\leq p\leq 4$
$$\delta(p):=\frac{1}{18}\left(1-\frac{2}{p}\right),\quad\varepsilon(p):=2\left(1-\frac{2}{p}\right).$$
\end{theo}
\begin{rema} Given a point $x_0$, we note that~\eqref{e:hyp-crit} and~\eqref{e:hyp-trans} are satisfied for an open and dense subset $\mathcal{U}_{x_0}$ of potentials in $\mathcal{C}^{\infty}(\IS^2,\IR)$ (endowed with its natural Fr\'echet topology). Assumption~\eqref{e:hyp-trans} implies that the Hamiltonian vector field is transverse to $\Gamma_{x_0}$ except at finitely many points. Combined with~\eqref{e:hyp-crit}, one has that, at the points where $X_{\langle V\rangle}(\gamma)$ is tangent to $\Gamma_{x_0}$, the tangency is of order $1$. See Remark~\ref{r:caustic} for an interpretation of these assumptions in terms of Lagrangian tori.
\end{rema}

\begin{rema} A direct Corollary of Theorem~\ref{t:maintheo} is that, if $K$ is a compact subset of $\IS^2$ such that, for every $x_0\in K$,~\eqref{e:hyp-crit} and~\eqref{e:hyp-trans} hold, then, for any solution $(\psi_\lambda,\lambda)$ to~\eqref{e:schrodinger},
\begin{equation}\label{e:Lp-compact}\|\psi_\lambda\|_{L^p(K)}\leq C_{K,p}\log(2+|\lambda|)^{\varepsilon(p)}(1+|\lambda|)^{\sigma_0(p)-\delta(p)}\|\psi_\lambda\|_{L^2(\mathbb{S}^2)}.\end{equation}
Yet, our main result does not allow to take $K=\IS^2$ as $\operatorname{Crit}(\mathcal{R}(V))$ cannot be empty.
\end{rema}

This Theorem yields a \emph{local} improvement for $p\neq 6$ over Sogge's upper bounds near certain points of $\IS^2$ which are \emph{independent} of the sequence $(\psi_{\lambda})_\lambda$ under consideration. The condition on these points are of purely dynamical nature and they depend on the subprincipal symbol of our operator. It may happen that Sogge's upper bounds are saturated for these operators but this can only occur away from points $x_0$ verifying~\eqref{e:hyp-crit} and~\eqref{e:hyp-trans}. The critical case $p_c=6$ could maybe be treated using similar ideas and the methods of Blair and Sogge to handle this exponent on nonpositively curved surfaces~\cite{So17, BlSo19}. Yet, this would probably require a much more delicate analysis than the one presented in this article.

Our hypothesis~\eqref{e:hyp-crit} and~\eqref{e:hyp-trans} are reminiscent from assumptions that appear when studying joint eigenfunctions of quantum completely integrable systems -- see~\cite[\S 1]{ToZe02} for a definition. For instance, the critical points involved in hypothesis~\eqref{e:hyp-crit} were used to obtain lower bounds by Toth in~\cite{To96} and by Toth-Zelditch in~\cite{ToZe02, ToZe03}. Similarly, assumption~\eqref{e:hyp-trans} was recently used by Galkowski--Toth~\cite{GaTo20} and by Tacy~\cite{Ta19b} to study the growth of $L^{\infty}$-norms of joint eigenfunctions. The main differences with these last references are that we handle every $p\neq 6$ and that we consider here eigenfunctions of the \emph{single} operator $-\Delta_{g_0}+V$. In fact recall from~\cite[Lemma~1]{Gu81} (see also~\cite{W77}) that there exists a unitary pseudodifferential operator $\mathcal{U}$ of order $0$ such that 
\begin{equation}\label{e:guillemin-weinstein}\mathcal{U}^{-1}\left(-\Delta_{g_0}+V\right)\mathcal{U}=-\Delta_{g_0}+V^\sharp,\end{equation}
where $[\Delta_{g_0},V^\sharp]=0$ and where the principal symbol of $V^\sharp$ is $\mathcal{R}(V)$. In other words, $-\Delta_{g_0}+V$ is the sum of two commuting pseudodifferential operators $\widehat{H}_1:=\mathcal{U}\Delta_{g_0}\mathcal{U}^{-1}$ and $\widehat{H}_2:=\mathcal{U}V^\sharp\mathcal{U}^{-1}$. In particular, it is a quantum completely integrable operator in the sense of~\cite[\S 1]{ToZe02} whenever $X_{\langle V\rangle}$ does not vanish on a dense and open subset of finite complexity (say outside finitely many points). Hence, upper bounds on $L^p$ norms of solutions to~\eqref{e:schrodinger} which are joint eigenfunctions of $(\widehat{H}_1,\widehat{H}_2)$ would follow from the results in~\cite{GaTo20, Ta19b} in the range $p>6$. However, in Theorem~\ref{t:maintheo}, we only suppose $p\neq 6$ and we do not make any assumption on the fact that $\psi_\lambda$ is a joint eigenfunction\footnote{We are not aware of a geometric criterion ensuring that all eigenfunctions are joint eigenfunctions. Yet, this is for instance achieved when the spectrum of $-\Delta_{g_0}+V$ is simple, which is the case for a residual set of potentials~\cite[Th.~7]{Uh76}.} of $(\widehat{H}_1,\widehat{H}_2)$ which makes the analysis slightly more delicate. Despite that, Theorem~\ref{t:maintheo} shows that there is room for (weaker) polynomial improvements on~\eqref{e:Lp-sogge} even for such eigenfunctions and even for $p<6$. In~\cite{Ta19}, Tacy obtained better estimates up to $p=2$ but she made stronger assumptions than ours on the sequence of eigenfunctions. Indeed, when restricted to our framework, the main result from this reference applies to sequences of joint eigenfunctions that concentrate away from the critical points of $\mathcal{R}(V)|_{\Gamma_{x_0}}$.

%Now, according to~\cite[Cor.~2.1.B]{Tay91}, $\|\mathcal{U}\|_{L^p\rightarrow L^p}<+\infty$ for every $1<p<\infty$ and  $\|\mathcal{U}\|_{\ml{C}^0\rightarrow \ml{C}^0}<+\infty$. 

\subsection{Earlier and related results}

The upper bounds~\eqref{e:Lp-sogge} are in fact valid in the general framework of compact Riemannian surfaces and, up to modifying the exponent $\sigma_0(p)$, they remain true in higher dimensions~\cite{So88}. Trying to improve them using the geometry of the manifold has been a classical topic in global harmonic analysis over the last thirty years.

\begin{itemize}
 \item \textbf{Flat tori.} In the case of flat tori and where $V\equiv 0$, this was achieved by Cooke~\cite{Co71} and Zygmund~\cite{Zy74} in dimension $2$ while the higher dimensional case was pursued by Bourgain~\cite{Bo93} and by Bourgain-Demeter~\cite{BoDe15}. In that case, one can use the arithmetic structure of the torus to get polynomial improvements over~\eqref{e:Lp-sogge}. See also~\cite{Wa11} for the case of Schr\"odinger operators on $2$-dimensional tori. To the best of the author's knowledge, flat tori is almost the only geometric framework where one can get global polynomial improvements without any further assumptions on the sequence of eigenfunctions (see below for the case of joint eigenfunctions). We can also mention~\cite{Zh20, Zh21} for recent improvements on compact Lie groups.
 \item \textbf{Negatively curved manifolds.} Another important class of examples where one expects improvements are negatively curved manifolds. For $p=\infty$, B\'erard showed how to get logarithmic improvements~\cite{Be77}. This logarithmic gain was extended to the range $p>p_c$ by Hassell and Tacy~\cite{HaTa15} and to manifolds without conjugate points by Bonthonneau~\cite{Bon17}. Still on negatively curved manifolds and for $p\leq p_c$, we obtained together with Hezari a logarithmic gain along generic sequences of eigenfunctions~\cite{HR16}. In a series of works related to Kakeya-Nikodym norms~\cite{So17, BlSo18, BlSo19}, Blair and Sogge proved logarithmic gains (with a slightly worst exponent) in this geometric context without any restriction on the sequence of eigenfunctions.
 \item \textbf{Arithmetic eigenfunctions.} A natural way to look for improvements over~\eqref{e:Lp-sogge} is to consider families of eigenfunctions that verify extra symmetries, for instance joint eigenfunctions of the Laplacian and of a family of commuting operators. In the case of a compact arithmetic surface, Iwaniec and Sarnak considered joint eigenfunctions of the Laplacian and of Hecke operators. For such sequences of eigenfunctions, they proved a polynomial improvement in the case of the $L^{\infty}$-norm~\cite{IwSar95}. In the case of the sphere, Brooks and Le Masson considered the related problem of joint eigenfunctions of $\Delta_{g_0}$ and the averaging operator for a finitely-generated free algebraic subgroup of $SO(3)$~\cite{BrLM20}. For such eigenfunctions, they obtained the same logarithmic improvement as Hassell and Tacy in the negatively curved case. On a rank $r$ symmetric space of dimension $n$, Sarnak improved the bound on the $L^{\infty}$-norm by a polynomial factor for eigenfunctions of the full ring of differential operators~\cite{Sar04}. This was generalized to the case of $L^p$-norms by Marshall~\cite{Ma16}.
 \item \textbf{Completely integrable systems.} Another context (closely related to ours) is the case of \emph{joint} eigenfunctions of a quantum completely integrable system. Toth and Zelditch proved that such eigenfunctions cannot have their $L^p$ norms uniformly bounded except in the case of flat tori~\cite{ToZe02, ToZe03}. See~\cite[Ch.~11]{Ze17} for a detailed discussion on joint eigenfunctions of quantum completely integrable systems. More recently, Galkowski and Toth obtained polynomial improvements on the $L^{\infty}$-bound for joint eigenfunctions of a quantum completely integrable systems~\cite{GaTo20} and Tacy proved improved Sogge's bounds for joint eigenfunctions of general families of semiclassical pseudodifferential operators~\cite{Ta19, Ta19b}.
 \item \textbf{Local improvements.} Sogge and Zelditch considered the problem from a more local perspective as we are doing here. They proved that, if, for a given point $x_0$ on a Riemannian manifold $(M,g)$, the set of covectors $\xi\in S_{x_0}^*M$ that come back to $x_0$ in finite time has zero measure, then one can improve locally near $x_0$ the upper bound on the $L^{\infty}$-norm by a $o(1)$ term~\cite{SoZe02}. This was based on improvements on the remainder in the local Weyl law. See also~\cite{Sa88} for earlier related results of Safarov. This result was later extended by Sogge, Toth and Zelditch under the weaker assumptions that the set of recurrent co-vectors at $x_0$ has $0$-measure\footnote{We emphasize that Theorem~\ref{t:maintheo} considers the somehow opposite case where the set of recurrent vectors has full measure. Despite that, we are able to get local polynomial improvements using the periodicity of the geodesic flow and the presence of a subprincipal symbol.}~\cite{SoToZe11}. We also refer to~\cite{SoZe16} for further developments of this approach when the metric is analytic and to~\cite[Ch.~10]{Ze17} for a detailed review. Related to these works, Galkowski and Toth showed how to relate precisely the growth of the $L^{\infty}$-norm near a point $x_0$ to the semiclassical measure restricted to the (geodesic) flow-out of the fiber $S_{x_0}^*M$~\cite{GaTo18} -- see also~\cite{Ga19}. More precisely, they proved that, if the $n$-dimensional Hausdorff measure of the support of this restriction is $0$, then one can get a $o(1)$-improvement on the growth of $L^{\infty}$-norm near $x_0$.
 \item \textbf{Using Gaussian beams.} This local approach was further improved by Canzani-Galkowski in a series of work using Gaussian beams~\cite{CaGa19, CaGa20}. In~\cite[Th.~1]{CaGa20}, they showed how to use this notion in order to give quantitative and at most logarithmic improvements on the growth of $L^p$-norms near a point $x_0$ when the conjugate points to $x_0$ do not pass too close to $x_0$. Among other things, they recover in that manner the results of B\'erard, Hassell-Tacy and Bonthonneau on manifolds without conjugate points. Besides that, they manage to deduce from their main results local improvements near $x_0$ on the growth of $L^p$ norms (for $p>p_c$) under quantitative assumptions on the geodesics passing through the point $x_0$ as in the works of Sogge, Toth and Zelditch. Finally, they also applied their main results to certain integrable (non-periodic) geometries on $\IS^2$ and obtain logarithmic improvements away from certain critical points when $p=\infty$~\cite[Th.~5]{CaGa19}. As in our framework, their result holds for the eigenfunctions of a single operator. 
\end{itemize}

\subsection{Strategy of proof}

In the range $p>6$, the proof is based on an argument to study the growth of $L^p$ norms that was used by Hezari and the author in~\cite{HR16} and further improved by Sogge in~\cite{So16}. It consists in relating the growth of $L^p$-norms to the growth of 
\begin{equation}\label{e:ave-ball}\int_{B_r(x)}|\psi_\lambda(y)|^2d\upsilon_{g_0}(y)\end{equation}
as $\lambda\rightarrow +\infty$ and $r\rightarrow 0^+$ (in a way that depends on $\lambda$). For $2\leq p<6$, we rather make use of results due to Blair and Sogge~\cite{So11, BlSo15, BlSo17} to control $L^p$-norms in terms of Kakeya-Nikodym averages around closed geodesics. See also~\cite{Bo09} for earlier related results of Bourgain. Then, we obtain rough bounds on these averages in terms of~\eqref{e:ave-ball}. The results from these references are briefly recalled (and adapted to Schr\"odinger eigenfunctions) in Sections~\ref{s:Lp} and~\ref{s:kakeya}. 

Up to smoothing the characteristic function of the balls, these local quantities can be interepreted in terms of Wigner distributions (or microlocal lifts). In particular, as was for instance observed by Shnirelman in his seminal work on quantum ergodicity~\cite{Sh74, Sh74b}, these distributions verify an almost invariance property by the geodesic flow. See for instance~\cite[Lemma 2, Eq.~(10)]{Sh74b}. This yields an upper bound of order $\mathcal{O}(r)$ on~\eqref{e:ave-ball} at least if $r$ does not go too fast to $0$ (say $r\gg\lambda^{-\frac{1}{2}}$). This is valid in a quite general framework. Yet, this is not sufficient to get an improvement over Sogge's upper bound. In order to implement this approach, one needs to have upper bounds of order $\mathcal{O}(r^{1+\alpha})$ for some $\alpha>0$, or at least $\mathcal{O}(\delta(r)r)$ with $\delta(r)\rightarrow 0$ as $r\rightarrow 0^+$. 

As pointed out by Sarnak in~\cite{Sar04}, a natural manner to look for improvements over Sogge's upper bounds is to consider operators commuting with the Laplacian and to study the $L^p$ norm of joint eigenfunctions. These joint eigenfunctions enjoy more symmetries which may lead to improvements. This was for instance the strategy followed in~\cite{IwSar95, BrLM20, Ma16, GaTo18, Ta19, Ta19b}. Here, we are not a priori in this situation as we consider eigenfunctions of the \emph{single} operator $-\Delta_{g_0}+V$ -- see the discussion following Theorem~\ref{t:maintheo}. However, the periodicity of the geodesic flow and the presence of the potential imply the existence of an extra invariance property besides the one by the geodesic flow. More precisely, in~\cite{MR16, MR19}, together with Maci\`a, we showed that Schr\"odinger eigenfunctions satisfy an extra invariance property by the Hamiltonian flow of $\mathcal{R}(V)$ which is reminiscent from the properties of joint eigenfunctions. This was achieved using Weinstein averaging method~\cite{W77}. Using this extra property, we will be able to get an upper bound of order $\mathcal{O}(r^{\frac{3}{2}})$ on~\eqref{e:ave-ball} up to scales $r\approx\lambda^{-\frac{2}{9}}$ near points verifying~\eqref{e:hyp-crit} and~\eqref{e:hyp-trans}. This will be the content of Section~\ref{s:invariance}. This additional invariance will be the reason for the polynomial improvement of Theorem~\ref{t:maintheo}. As we shall see in our proof\footnote{See for instance~\eqref{e:range-p}.}, the reason for being limited to $p\neq 6$ comes from this exponent $3/2$ and, in dimension $2$, any bound on~\eqref{e:ave-ball} of order $\mathcal{O}(r^{1+\alpha})$ with $\alpha>1/2$ would give a local improvement over Sogge's upper bound~\eqref{e:Lp-sogge} even for $p=6$ (using the arguments of \S~\ref{s:Lp}). 

%\begin{rema} Even if it concerns the eigenfunctions of the single operator $-\Delta_{g_0}+V$, this extra invariance property that we obtained with Maci\`a is reminiscent from the properties of joint eigenfunctions for completely integrable systems. The careful reader can verify that our proof would work as well for joint eigenfunctions of $-\Delta_{g_0}$ and
%$$\frac{1}{2\pi}\int_{0}^{2\pi}e^{-is \sqrt{-\Delta_{g_0}+\frac{1}{4}}}Ve^{is \sqrt{-\Delta_{g_0}+\frac{1}{4}}}ds.$$
%Yet, the bounds that we would obtain would not be better that the ones for joint eigenfunctions obtained by Tacy in~\cite{Ta19, Ta19b} or by Galkowski and Toth in~\cite{GaTo20}. Recall that, except on flat tori, they cannot be made uniformly bounded as shown by Toth and Zelditch~\cite{To96, ToZe02, ToZe03}. 
%\end{rema}

\section*{Acknowledgements}
I would like to address my warmest thanks to Hamid Hezari and Fabricio Maci\`a for my joint works with them~\cite{HR16, MR16, MR19} and for their many insights on these topics. I also thank Xiaolong Han for pointing me reference~\cite{Uh76} regarding the generic simplicity of the spectrum of Schr\"odinger operators and the anonymous referee for helpful suggestions. This work was supported by the Institut Universitaire de France and by the Agence Nationale de la Recherche through the PRC projects ODA (ANR-18-CE40-0020) and ADYCT (ANR-20-CE40-0017).

\section{Reduction to $L^2$ localized estimates for $p>6$}
\label{s:Lp}

In this section, we revisit an argument due to Sogge\footnote{See also~\cite{HR16} for earlier related arguments of Hezari and the author using semiclassical methods~\cite[\S 10]{Zw12}.} in order to relate $L^p$ estimates to localized $L^2$-estimates in small balls. This argument will allow us to get our upper bounds in the range $6<p\leq\infty$. The proof given in~\cite{So16} was for Laplace eigenfunctions and we verify that it can be adapted to Schr\"odinger eigenfunctions. 

\begin{rema} Due to our $L^2$-localized estimates in Section~\ref{s:invariance}, we could as well work only with $p=\infty$ and conclude by interpolation with the case $p=6$ in~\eqref{e:Lp-sogge}. Yet, we write things down for general $p$ in order to identify the quantitative improvements one would need to reach the case $p=6$. See Equation~\eqref{e:range-p} below.
\end{rema}

Let $\psi_\lambda$ be a solution to~\eqref{e:schrodinger} that we suppose to be $L^2$-normalized. In the following, we suppose that $\lambda^2$ is large enough so that we can pick $\lambda>0$. Following~\cite[\S 2]{So16} and for $j\in\mathbb{Z}_+$, we denote by $E_j$ the spectral projector onto the eigenspace of $\sqrt{-\Delta_{g_0}}$ with eigenvalue $\lambda_j:=\sqrt{j(j+1)}$. We write
\begin{equation}\label{e:quasimode-explicit}(\sqrt{-\Delta_{g_0}}-\lambda)\psi_\lambda=-(\sqrt{-\Delta_{g_0}}+\lambda)^{-1}V\psi_\lambda.\end{equation}
In particular, one has
\begin{equation}\label{e:quasimode}\left\|(\sqrt{-\Delta_{g_0}}-\lambda)\psi_\lambda\right\|_{L^2}=\left(\sum_{j\in\IZ^+}\frac{1}{(\lambda+\lambda_j)^2}\|E_j(V\psi_\lambda)\|^2\right)^{\frac{1}{2}}\leq \frac{1}{\lambda}\|V\psi_{\lambda}\|_{L^2}\leq \frac{\|V\|_{L^\infty}}{\lambda}.\end{equation}

%and we set
%$$\chi_{\lambda}:=\sum_{\lambda_j\in[\lambda,\lambda+1)}E_j.$$
%This is the projector on unit bands of frequencies. 
We also fix a nonnegative $\rho\in\mathcal{S}(\IR)$ satisfying
\begin{equation}\label{e:rho}
\rho(0)=1\quad\text{and}\quad\text{supp}(\hat{\rho})\subset[-1,1],
\end{equation}
where $\hat{\rho}$ is the Fourier transform of $\rho$. For $\lambda>0$ and $0<r\leq 1$, setting
$$T_{\lambda,r}:=\frac{1}{\pi}\int_{-\infty}^{+\infty}r^{-1}\hat{\rho}(r^{-1}t)e^{it\lambda}\cos(t\sqrt{-\Delta_{g_0}})dt,$$
one finds
$$T_{\lambda,r}=\rho\left(r\left(\lambda-\sqrt{-\Delta_{g_0}}\right)\right)+\rho\left(r\left(\lambda+\sqrt{-\Delta_{g_0}}\right)\right).$$
The main result of~\cite[Eq.~(3.1)]{So16} is that, for every $p>2$ and for every $f\in L^2(\mathbb{S}^2)$,
\begin{equation}\label{e:sogge}
 \|T_{\lambda,r}f\|_{L^p(\IS^2)}\leq C_p r^{-\frac{1}{2}}\lambda^{\sigma_0(p)}\|f\|_{L^2(\IS^2)},\quad\lambda\geq 1,\quad\lambda^{-1}\leq r\leq\frac{\pi}{2},
\end{equation}
where the constant $C_p$ is uniform for $(\lambda,r)$ in the above range. Recall now from Huygens principle that the Schwartz kernel $\cos(t\sqrt{-\Delta_{g_0}})(x,y)$ vanishes if the geodesic distance between $x$ and $y$ is $>t$. In particular, the Shwartz kernel $T_{\lambda,r}(x,y)$ of $T_{\lambda,r}$ vanishes if $d_{g_0}(x,y)>r$ thanks to our assumptions~\eqref{e:rho} on the support of $\rho$. Gathering these informations, Sogge observed that, for every $p>2$ and for every $f\in L^2(\mathbb{S}^2)$,
\begin{equation}\label{e:sogge2}
 \|T_{\lambda,r}f\|_{L^p(B_r(x_0))}\leq C_p r^{-\frac{1}{2}}\lambda^{\sigma_0(p)}\|f\|_{L^2(B_{2r}(x_0))},\quad\lambda\geq 1,\quad\lambda^{-1}\leq r\leq\frac{\pi}{2},
\end{equation}
where the constant $C_p$ is uniform for $(\lambda,r)$ in the above range and for $x_0\in\IS^2$. This will be referred as the Sogge's local $L^p$-estimate. Fix now some compact subset $K$ of $\mathbb{S}^2$. We can cover $K$ by finitely many balls $(B_r(x_l))_{l=1,\ldots, N(r)}$ of radius $r$ and centered at points inside $K$. We require that the number $N(r)$ is of order $\sim r^{-2}$ and that each point of $K$ is contained in at most $C_0$ balls of the covering $(B_{2r}(x_l))_{l=1,\ldots, N(r)}$. Here $C_0>0$ is independent of $r$ -- see for instance~\cite[Lemma~2]{CM11}. Recall that we have in mind to apply this result when $K=B_{r_0}(x_0)$ is a fixed ball. Hence, one has, for $2<p<\infty$ and for $f$ in $L^2(\IS^2)$,
\begin{eqnarray*}
 \|f\|_{L^p(K)}^p %&\leq & 2^{p-1}\left(\|T_{\lambda,r}f\|_{L^p(K)}^p+\|(T_{\lambda,r}-\text{Id})f\|_{L^p(\IS^2)}^p\right)\\
 &\leq & 2^{p-1}\left(\sum_{l=1}^{N(r)}\|T_{\lambda,r}f\|_{L^p(B_r(x_l))}^p+\|(T_{\lambda,r}-\text{Id})f\|_{L^p(\IS^2)}^p\right)\\
  &\leq & C_p r^{-\frac{p}{2}}\lambda^{\sigma_0(p)p}\sum_{l=1}^{N(r)}\|f\|_{L^2(B_{2r}(x_l))}^p+C_p\|(T_{\lambda,r}-\text{Id})f\|_{L^p(\IS^2)}^p\\
   %&\leq & C_p r^{-\frac{p}{2}}\lambda^{\sigma_0(p)p}\left(\max_{1\leq l\leq N(r)}\left\{\|f\|_{L^2(B_{2r}(x_l))}^{p-2}\right\}\right)\sum_{l=1}^{N(r)}\|f\|_{L^2(B_{2r}(x_l))}^2\\
   %&+&C_p\|(T_{\lambda,r}-\text{Id})f\|_{L^p(\IS^2)}^p\\
    &\leq & C_pC_0 r^{-\frac{p}{2}}\lambda^{\sigma_0(p)p}\left(\max_{1\leq l\leq N(r)}\left\{\|f\|_{L^2(B_{2r}(x_l))}^{p-2}\right\}\right)\|f\|_{L^2(\IS^2)}^2   +C_p\|(T_{\lambda,r}-\text{Id})f\|_{L^p(\IS^2)}^p.
\end{eqnarray*}
Hence, one finds
\begin{lemm}\label{l:step-lemma} Let $K$ be a compact subset of $\IS^2$ and let $(B_r(x_l))_{l=1,\ldots, N(r)}$ be a cover of $K$ with the above properties. Then, one has 
\begin{equation}\label{e:lpbound}
 \|f\|_{L^p(K)}\leq C_p'\left(r^{-\frac{1}{2}}\lambda^{\sigma_0(p)}\left(\max_{1\leq l\leq N(r)}\left\{\|f\|_{L^2(B_{2r}(x_l))}^{1-\frac{2}{p}}\right\}\right)\|f\|_{L^2(\IS^2)}^2
   +\|(T_{\lambda,r}-\text{Id})f\|_{L^p(\IS^2)}\right).
   \end{equation}\end{lemm}
This upper bound is valid uniformly in the range $\lambda\geq 1$ and $ \lambda^{-1}\leq r\leq\frac{\pi}{2}.$ Similarly, in the case of the $L^{\infty}$ norm, we would get
\begin{equation}\label{e:linftybound}\|f\|_{L^\infty(K)}\leq Cr^{-\frac{1}{2}}\lambda^{\frac{1}{2}}\left(\max_{1\leq l\leq N(r)}\left\{\|f\|_{L^2(B_{2r}(x_l))}\right\}\right)+\|(T_{\lambda,r}-\text{Id})f\|_{L^\infty(\IS^2)}.
\end{equation}
Note that so far we did not use the eigenvalue equation~\eqref{e:quasimode-explicit} and this is valid for any $f$ in $L^2(\IS^2)$. We will now specify these results in the case where $f=\psi_\lambda$. We begin with the remainder term:
\begin{prop}\label{p:remainder} Let $2<p\leq\infty$ and let $0< \beta<1$. Then, there exists a constant $C>0$ such that, for any solution $\psi_\lambda$ to~\eqref{e:schrodinger} with $\lambda\geq 1$ and for any $\lambda^{-\beta}\leq r\leq\frac{\pi}{2}$, one has
 $$\left\|\left(T_{\lambda,r}-\text{Id}\right)\psi_\lambda\right\|_{L^p(\mathbb{S}^2)}\leq C(r\lambda)^{\sigma_0(p)} \|\psi_\lambda\|_{L^2(\mathbb{S}^2)},$$
\end{prop}
Regarding Lemma~\ref{l:step-lemma} which already used the Sogge's local $L^p$-estimate, this proposition is the additional ingredient we need to take into account the terms coming from the potential $V$. Gathering this Proposition with our estimates~\eqref{e:lpbound} and~\eqref{e:linftybound} on $\|f\|_{L^p(K)}$, we find that, for $\lambda\geq 1$, $\lambda^{-\beta}\leq r\leq \frac{\pi}{2}$ (with $\beta<1$), for any $2< p\leq +\infty$ and for any $L^2$-normalized solution $\psi_\lambda$ to~\eqref{e:schrodinger},
\begin{equation}\label{e:localized-Lp}
\|\psi_{\lambda}\|_{L^p(K)}\leq C_p \left(r^{-\frac{1}{2}}\lambda^{\sigma_0(p)}\max_{1\leq l\leq N(r)}\left\{\|\psi_\lambda\|_{L^2(B_{2r}(x_l))}^{1-\frac{2}{p}}\right\}+(r\lambda)^{\sigma_0(p)}\right).
\end{equation}
The involved constants $C_p>0$ depend only on $V$, $K$, $\rho$, $\beta$ and $p$. Hence, as in~\cite{HR16, So16}, we have reduced the problem of estimating the $L^p$ norm of Schr\"odinger eigenfunctions to determining bounds on $L^2$-localized norms,
\begin{equation}\label{e:smallball}\int_{B_{2r}(x_l)}|\psi_\lambda(x)|^2d\upsilon_{g_0}(x),\end{equation}
as $\lambda\rightarrow +\infty$ with $r$ verifying $\lambda^{-\beta}\leq r\leq \frac{\pi}{2}$. In particular, if, for some $0<\alpha\leq 1$, we were able to bound~\eqref{e:smallball} uniformly (in terms of $\lambda$) by $Cr^{1+\alpha}$, then we would be able to get an improved upper bound inside $K$ of the form
$$\|\psi_{\lambda}\|_{L^p(K)}\leq C_{p,K}\left(r^{\frac{\alpha}{2}-\frac{1+\alpha}{p}}\lambda^{\sigma_0(p)}+ (r\lambda)^{\sigma_0(p)}\right),$$
in the range 
\begin{equation}\label{e:range-p}\frac{\alpha}{2}-\frac{1+\alpha}{p}>0\quad\Longleftrightarrow\quad p>2\left(1+\frac{1}{\alpha}\right).\end{equation}
%Hence, if we could prove this uniform control for radii $r_\lambda\geq \lambda^{-\beta}$ that tend to $0$ as $\lambda\rightarrow+\infty,$ then we would get an improvement on the $L^p$ norm. 
%This was exactly the strategy that we followed with Hezari in the framework of negatively curved surfaces~\cite{HR16} where we were able to reach $\alpha=1$ and $r_\lambda=|\log\lambda|^{-\frac{1}{4}+o(1)}$ for a density $1$ subsequence of eigenfunctions and for $K=M$. 
However, as explained in~\cite[\S 4]{So16}, one cannot expect such improved bounds on the sphere when $V\equiv 0$ thanks to the example of the spherical harmonics. In section~\ref{s:invariance}, we shall see how to get \emph{locally} improved bounds on~\eqref{e:smallball} when $V$ does not identically vanish. Before going to this question, we give the proof of Proposition~\ref{p:remainder}.
\begin{proof}The ideas of the proof are standard (see e.g.~\cite{So16}) and we detail them for the sake of completeness.
Considering a solution to~\eqref{e:quasimode-explicit} and letting $2\leq p\leq +\infty$, one has
\begin{eqnarray*}\left\|\left(T_{\lambda,r}-\text{Id}\right)\psi_\lambda\right\|_{L^p(\mathbb{S}^2)}&\leq&\sum_{j\in\mathbb{Z}_+}\left\|E_j\left(T_{\lambda,r}-\text{Id}\right)E_j\psi_\lambda\right\|_{L^p(\mathbb{S}^2)}\\
% &\leq &\sum_{j\in\mathbb{Z}_+}\left|\rho(r(\lambda-\lambda_j))-1+\rho(r(\lambda+\lambda_j))\right|\left\|E_j(\psi_\lambda)\right\|_{L^p(\mathbb{S}^2)}\\
 &\leq &\sum_{j\in\mathbb{Z}_+}\left(\left|\rho(r(\lambda-\lambda_j))-1\right|+\left|\rho(r(\lambda+\lambda_j))\right|\right)\left\|E_j(\psi_\lambda)\right\|_{L^p(\mathbb{S}^2)}.
\end{eqnarray*}
As $\rho$ belongs to the Schwartz class, we find using Sogge's estimate~\eqref{e:Lp-sogge} that, for every $N\geq 1$, there exists $C_N>0$ such that, for $\lambda\geq 1$ and $r\geq\lambda^{-\beta}$,
$$\sum_{j\in\mathbb{Z}_+}\left|\rho(r(\lambda+\lambda_j))\right|\left\|E_j(\psi_\lambda)\right\|_{L^p(\mathbb{S}^2)}\leq C_N(1+r\lambda)^{-N}\|\psi_\lambda\|_{L^2(\mathbb{S}^2)}.$$
%Hence, one gets
%$$\left\|\left(T_{\lambda,r}-\text{Id}\right)\psi_\lambda\right\|_{L^p(\mathbb{S}^2)}\leq\sum_{j\in\mathbb{Z}_+}\left|\rho(r(\lambda-\lambda_j))-1\right|\left\|E_j(\psi_\lambda)\right\|_{L^p(\mathbb{S}^2)}+C_N(1+r\lambda)^{-N}\|\psi_\lambda\|_{L^2(\mathbb{S}^2)}.$$
Using one more time Sogge's estimate, we deduce that
\begin{equation}\label{e:remainder}\left\|\left(T_{\lambda,r}-\text{Id}\right)\psi_\lambda\right\|_{L^p(\mathbb{S}^2)}\leq\sum_{j\in\mathbb{Z}_+}\left|\rho(r(\lambda-\lambda_j))-1\right|\lambda_j^{\sigma_0(p)}\left\|E_j(\psi_\lambda)\right\|_{L^2(\mathbb{S}^2)}+C_N(1+r\lambda)^{-N}\|\psi_\lambda\|_{L^2(\mathbb{S}^2)}.\end{equation}
We now fix some $\delta\geq r$ so that $\delta\leq r\lambda$ and we split the sum over $j\in\mathbb{Z}_+$ in two parts. On the one hand, we consider the $j$ such that $|\lambda-\lambda_j|\leq \delta/r$ and on the other hand, the integers such that $|\lambda-\lambda_j|> \delta/r$. Recall that $\lambda_j^2=j(j+1)$. Hence, the number of terms in the first sum is $\mathcal{O}(\delta/r)$ and one is left with
\begin{eqnarray*}
 \left\|\left(T_{\lambda,r}-\text{Id}\right)\psi_\lambda\right\|_{L^p(\mathbb{S}^2)}&\leq&\sum_{j\in\mathbb{Z}_+: |\lambda-\lambda_j|> \delta/r}\left|\rho(r(\lambda-\lambda_j))-1\right|\lambda_j^{\sigma_0(p)}\left\|E_j(\psi_\lambda)\right\|_{L^2(\mathbb{S}^2)}\\
 &+& \left(C\frac{\delta^2}{r}\lambda^{\sigma_0(p)}+C_N(1+r\lambda)^{-N}\right)\|\psi_\lambda\|_{L^2(\mathbb{S}^2)}.
\end{eqnarray*}
For the remaining sum, we can finally make use of the eigenvalue equation~\eqref{e:quasimode-explicit}. It implies the existence of some constant $C_{\rho,V}>0$ depending only on $\rho$ and $V$ such that
$$\sum_{j\in\mathbb{Z}_+:|\lambda-\lambda_j|> \delta/r}\left|\rho(r(\lambda-\lambda_j))-1\right|\lambda_j^{\sigma_0(p)}\left\|E_j(\psi_\lambda)\right\|_{L^2(\mathbb{S}^2)}\leq C_{\rho,V}\sum_{j\in\mathbb{Z}_+:|\lambda-\lambda_j|> \delta/r}\frac{\lambda_j^{\sigma_0(p)}}{|\lambda^2-\lambda_j^2|}\|\psi_{\lambda}\|_{L^2(\IS^2)}.$$
As $\sigma_0(p)$ varies between $0$ (for $p=2$) and $1/2$ (for $p=\infty$), this last quantity is finite and it remains to evaluate
\begin{equation}\label{e:bad-remainder}\sum_{j\in\mathbb{Z}_+:|\lambda-\lambda_j|> \delta/r}\frac{\lambda_j^{\sigma_0(p)}}{|\lambda^2-\lambda_j^2|}\end{equation}
in terms of $\delta$, $r$, $\lambda$ and $p$. We now recall that, for $X>0$, one has $(1+X)^{\sigma_0(p)}\leq 1+X^{\sigma_0(p)}$ (as $\sigma_0(p)\leq1/2$). Hence, one has
\begin{eqnarray*}\sum_{j\in\mathbb{Z}_+:|\lambda-\lambda_j|> \delta/r}\frac{\lambda_j^{\sigma_0(p)}}{|\lambda^2-\lambda_j^2|}&\leq& \sum_{j\in\mathbb{Z}_+:|\lambda-\lambda_j|> \delta/r}\frac{|\lambda-\lambda_j|^{\sigma_0(p)}}{|\lambda^2-\lambda_j^2|}+\sum_{j\in\mathbb{Z}_+:|\lambda-\lambda_j|> \delta/r}\frac{\lambda^{\sigma_0(p)}}{|\lambda^2-\lambda_j^2|}\\
 &\leq& 2\sum_{j\in\mathbb{Z}_+:|\lambda-\lambda_j|> \delta/r}\frac{\lambda^{-1+\frac{3}{2}\sigma_0(p)}}{|\lambda-\lambda_j|^{1+\frac{\sigma_0(p)}{2}}}\\
 &\leq & 2\lambda^{-\frac{1}{4}}\sum_{j\in\mathbb{Z}_+:|\lambda-\sqrt{j(j+1)}|> \delta/r}\frac{1}{|\lambda-\sqrt{j(j+1)}|^{1+\frac{\sigma_0(p)}{2}}}\\
 %& \leq &2\lambda^{-\frac{1}{4}}\sum_{j\in\mathbb{Z}_+:j+1<\lambda-1}\frac{1}{(\lambda-(j+1))^{1+\frac{\sigma_0(p)}{2}}}\\
 %&+&2\lambda^{-\frac{1}{4}}\sum_{j\in\mathbb{Z}_+:j>\lambda+1}\frac{1}{(j-\lambda)^{1+\frac{\sigma_0(p)}{2}}}\\
 &\leq&C\lambda^{-\frac{1}{4}}\sum_{j\in\mathbb{Z}_+^*}j^{-1-\frac{\sigma_0(p)}{2}}.
\end{eqnarray*}
In summary, if we suppose that $r\geq\lambda^{-\beta}$ (for some $\beta<1$), we obtain the following upper bound
$$\left\|\left(T_{\lambda,r}-\text{Id}\right)\psi_\lambda\right\|_{L^p(\mathbb{S}^2)}\leq C\left(\frac{\delta^2}{r}\lambda^{\sigma_0(p)}+\lambda^{-\frac{1}{4}}\right) \|\psi_\lambda\|_{L^2(\mathbb{S}^2)},$$
where $C>0$ depends on $\rho$, $V$, $\beta$ and $p$. Recall that we supposed $r\leq \delta\leq r\lambda$. Hence, as $0\leq\sigma(p)\leq \frac{1}{2}$, we can set $\delta=r^{\frac{1+\sigma_0(p)}{2}}$ provided $r\geq \lambda^{-\frac{2}{\sigma_0(p)+1}}$, which is ensured by our assumption $r\geq\lambda^{-\beta}$. Implementing this, we obtain the existence of a constant $C_{\rho,V,\beta,p}>0$ such that
$$\left\|\left(T_{\lambda,r}-\text{Id}\right)\psi_\lambda\right\|_{L^p(\mathbb{S}^2)}\leq C_{\rho,V,\beta,p}(r\lambda)^{\sigma_0(p)} \|\psi_\lambda\|_{L^2(\mathbb{S}^2)},$$
as long as $r\geq \lambda^{-\beta}$.
\end{proof}

\begin{rema}\label{r:Lp-semiclassical} In view of applications of our method to semiclassical problems, it is worth noting that the above arguments work as well for solutions to
\begin{equation}\label{e:schrodinger-pert-lambda}-\Delta_{g_0}\psi_\lambda+\beta_\lambda V\psi_\lambda =\lambda^2\psi_\lambda,\quad\|\psi_{\lambda}\|_{L^2(\IS^2)}=1,\end{equation}
 where $(\beta_\lambda)_{\lambda}$ is a given nonnegative sequence that may tend to $+\infty$. In that case, the upper bound~\eqref{e:localized-Lp} becomes, for every $\epsilon>0$,
 \begin{equation}\label{e:localized-Lp-bis}
\|\psi_{\lambda}\|_{L^p(K)}\leq C_{p,\epsilon} \left(r^{-\frac{1}{2}}\lambda^{\sigma_0(p)}\max_{1\leq l\leq N(r)}\left\{\|\psi_\lambda\|_{L^2(B_{2r}(x_l))}^{1-\frac{2}{p}}\right\}+(r\lambda)^{\sigma_0(p)}+\beta_\lambda \lambda^{-1+\epsilon}\lambda^{\sigma_0(p)}\right).
\end{equation}
The calculation is indeed exactly the same except for the upper bound on the size of the remainder in~\eqref{e:bad-remainder} that we need to improve. Hence, we have potentially improvements as long as\footnote{This can probably sligthly improved to replace the $\lambda^{\epsilon}$ by some logarithmic factor but we did not try to optimize that.} $\beta_\lambda \lambda^{-1+\epsilon}\rightarrow 0$.
\end{rema}

\section{Reduction to $L^2$ localized estimates for $p<6$ via Kakeya-Nikodym bounds}
\label{s:kakeya}

We now deal with the range $2<p<6$ which can also be reduced to estimating similar quantities. For such $p$, we can make use of the results of Blair and Sogge relating the growth of $L^p$ norms for small $p$ to Kakeya-Nikodym averages.

We let $0\leq \chi\leq 1$ be a smooth cutoff function which is equal to $1$ on $[-1,1]$ and to $0$ outside $[-2,2]$. Given $x\in \IS^2$, we denote by $\exp_{x}$ the exponential map induced by the metric $g_0$ and we set
$$\chi_{x,r}(y):=\chi\left(\frac{\|\exp_{x}^{-1}(y)\|}{r}\right)\in\mathcal{C}^{\infty}(\IS^2).$$
This function is equal to $1$ on $B_r(x)$ and to $0$ outside $B_{2r}(x)$. We fix some $r_0>0$ and some $x_0\in\IS^2$. For any \emph{normalized} solution to~\eqref{e:schrodinger}, one has 
$$-\lambda^{-2}\Delta_{g_0}\psi_\lambda-\psi_{\lambda}=\lambda^{-2}V\psi_{\lambda}.$$
In particular, one can verify, using commutation rules for semiclassical pseudodifferential operators~\cite[\S~4 and~14]{Zw12},
\begin{equation}\label{e:KN-quasimode}(-\lambda^{-2}\Delta_{g_0}-1)^k\left(\chi_{x_0,r_0}\psi_\lambda\right)=\mathcal{O}(\lambda^{-k}),\quad k=1,2.\end{equation}
These two assumptions are exactly the ones needed to apply~\cite[Th.~1.1]{BlSo17} in dimension $2$. In order to formulate this result, we denote by $\tilde{G}(\IS^2)$ the set of unit length geodesic segments in $\IS^2$ and, for every $r>0$ and for every $\gamma\in \tilde{G}(\IS^2)$,
$$\mathcal{T}_r(\gamma):=\left\{x\in\IS^2:\ d_{g_0}(x,\gamma)\leq r\right\}.$$
With these conventions, the main result from~\cite{BlSo17} applied to $\chi_{x_0,r_0}\psi_\lambda$ tells us that, for $4<p<6$,
\begin{equation}\label{e:KN-no-log}\left\|\psi_{\lambda}\right\|_{L^p(B_{r_0}(x_0))}\leq C_p\lambda^{\sigma_0(p)}\left(\sup_{\gamma\in \tilde{G}(\IS^2)}\int_{B_{2r_0}(x_0)\cap \mathcal{T}_{\lambda^{-\frac{1}{2}}}(\gamma)}|\psi_\lambda(x)|^2d\upsilon_{g_0}(x)\right)^{\frac{1}{2}\left(\frac{6}{p}-1\right)},\end{equation}
and
\begin{equation}\label{e:KN-log}\left\|\psi_{\lambda}\right\|_{L^4(B_{r_0}(x_0))}\leq C_p(\log\lambda)\lambda^{\frac{1}{8}}\left(\sup_{\gamma\in \tilde{G}(\IS^2)}\int_{B_{2r_0}(x_0)\cap \mathcal{T}_{\lambda^{-\frac{1}{2}}}(\gamma)}|\psi_\lambda(x)|^2d\upsilon_{g_0}(x)\right)^{\frac{1}{4}},\end{equation}
where the constants $C_p>0$ depend only on $p$. These kinds of upper bounds are referred to as Kakeya-Nikodym bounds. They were initially introduced by Bourgain~\cite{Bo09} and further developped by Sogge~\cite{So11, So17} and Blair-Sogge~\cite{BlSo15, BlSo17, BlSo18, BlSo19}. One of the main objectives is to reduce (at least for small $p$) improvements on $L^p$-estimates to $L^2$-estimates on tubular neighborhoods of geodesics. This strategy culminated in~\cite{BlSo19} where logarithmic improvements on Sogge's $L^p$ estimates were obtained on nonpositively curved manifolds for every $p\leq p_c$. As we shall see below, this strategy remains efficient for integrable geometries where we can also analyze the $L^2$-mass near geodesics via averaging methods.

Thanks to these results, it is sufficient to derive nontrivial upper bounds on the Kakeya-Nikodym averages
$$\int_{B_{2r_0}(x_0)\cap \mathcal{T}_{\lambda^{-\frac{1}{2}}}(\gamma)}|\psi_\lambda(x)|^2d\upsilon_{g_0}(x).$$
in order to improve locally Sogge's upper bounds~\eqref{e:sogge} in the range $4<p<6$. By interpolation, it will automatically yields an improvement for $2<p<4$. 

Finally, we can relate these quantities to the ones appearing in~\eqref{e:smallball}. Indeed, we can pick $0<\beta<1/2$ and we can cover $B_{2r_0}(x_0)\cap \mathcal{T}_{\lambda^{-\frac{1}{2}}}(\gamma)$ by a family of $2r_0r^{-1}$ balls of radius $r\geq \lambda^{-\beta}$ centered on a point of $\gamma\cap B_{2r_0}(x_0)$. Hence, one has 
\begin{equation}\label{e:KN-smallball}
 \int_{B_{2r_0}(x_0)\cap \mathcal{T}_{\lambda^{-\frac{1}{2}}}(\gamma)}|\psi_\lambda(x)|^2d\upsilon_{g_0}(x)\leq 4r_0r^{-1}\sup_{x\in\gamma\cap B_{2r_0}(x_0)}\left\{\int_{B_{r}(x)}|\psi_\lambda(y)|^2d\upsilon_{g_0}(y)\right\},
\end{equation}
which are exactly the quantities that appeared in Section~\ref{s:Lp}. Hence, in both cases, we are reduced to estimating these localized $L^2$-estimates.

\begin{rema}\label{r:KN-Lp}
 As in Remark~\ref{r:Lp-semiclassical}, we can consider solutions to~\eqref{e:schrodinger-pert-lambda}. One can verify that the assumption~\eqref{e:KN-quasimode} is still verified as long as $0\leq \beta_\lambda\leq\lambda$. Hence,~\eqref{e:KN-no-log} and~\eqref{e:KN-log} remain true in that generalized framework.
\end{rema}

\begin{rema} As we will only consider balls of radius $r\gg\lambda^{-\frac{1}{2}}$, the logarithmic factor appearing in~\eqref{e:KN-log} could probably be removed following~\cite{BlSo15}. 
\end{rema}

\section{$L^2$-localized estimates using invariance by the classical flows}\label{s:invariance}

Thanks to~\eqref{e:localized-Lp},~\eqref{e:KN-no-log},~\eqref{e:KN-log} and~\eqref{e:KN-smallball}, we know that proving Theorem~\ref{t:maintheo} amounts to control uniformly the following quantity
$$M_{B_{r_0}(x_0), \alpha,r}(\psi_\lambda):=\sup\left\{\frac{1}{r^{1+\alpha}}\int_{B_{r}(x)}|\psi_\lambda(y)|^2d\upsilon_{g_0}(y):x\in B_{r_0}(x_0)\right\},$$
with $0<\alpha\leq 1$ and $\lambda^{-\beta}\leq r$ that goes to $0$ as $\lambda\rightarrow +\infty$. The following Proposition answers this problem and it is the main new technical result of the article:
\begin{prop}\label{p:main-prop} Let $x_0$ be a point in $\IS^2$ verifying the assumption of Theorem~\ref{t:maintheo}. Then, there exist $r_0>0$ and $C_0>0$ such that, for any $(\psi_\lambda,\lambda)$ solution to~\eqref{e:schrodinger},
 $$\lambda^{-\frac{2}{9}}\leq r\leq \frac{\pi}{2}\quad\Longrightarrow\quad M_{B_{r_0}(x_0),\frac{1}{2},r}\left(\psi_\lambda\right)\leq C_0 \|\psi_\lambda\|_{L^2(\IS^2)}^2.$$
\end{prop}

\begin{rema}
The exponent in $\lambda^{-2/9}$ appears as follows in the argument. On the one hand, we use semiclassical arguments for exotic class of symbols (with $\lambda^{-\beta}$ loss in the derivatives) and this yields remainder terms of size $\ml{O}(\lambda^{-1+3\beta})$. This semiclassical part of the argument is based on Egorov and composition theorems and the remainders cannot be drastically improved. See for instance~\eqref{e:final-bound-ball}. On the other hand, we need to estimate classical averages by some Hamiltonian flow and this is where we use in an essential way our assumption on the potential $V$. Without these assumptions, we would get a crude upper bound $\ml{O}(r)$ and these hypothesis allow to upgrade this bound to $\ml{O}(r^{1+\frac{1}{2}})$. See for instance~\eqref{e:bound-tangent-case}. Here the fact that $\ml{R}(V)|_{\Gamma_{x_0}}$ is a Morse function implies that the tangency have order at most $1$. For higher order tangencies, we would have probably obtained some slightly worst bound $\ml{O}(r^{1+\frac{1}{k}})$ (for some large enough $k$) at the expense of some extra tedious work. In the end, we take $r$ such that $\ml{O}(r^{1+\frac{1}{2}})$ and $\ml{O}(\lambda^{-1+3\beta})$ are of the same order which yields the exponent $2/9$.
\end{rema}

Implementing this bound in~\eqref{e:localized-Lp} and in~\eqref{e:KN-no-log}, we find that, for $4<p\leq\infty$ and for $\lambda>0$,
$$\left\|\psi_{\lambda}\right\|_{L^p(B_{r_0}(x_0))}\leq C_{p,x_0}\lambda^{\sigma_0(p)-\frac{1}{18}\left|1-\frac{6}{p}\right|}\|\psi_\lambda\|_{L^2(\IS^2)}.$$
%Using~\eqref{e:KN-no-log}, we get on the range $4<p<6$ and for $\lambda>0$,
%$$\left\|\psi_{\lambda}\right\|_{L^p(B_{r_0}(x_0))}\leq C_{p,x_0}\lambda^{\sigma_0(p)-\frac{1}{18}\left(\frac{6}{p}-1\right)}\|\psi_\lambda\|_{L^2(\IS^2)}.$$
Finally, for $p=4$, we derive from~\eqref{e:KN-log} that, for $\lambda>1$,
$$\left\|\psi_{\lambda}\right\|_{L^4(B_{r_0}(x_0))}\leq C_{4,x_0}(\log\lambda) \lambda^{\frac{1}{8}-\frac{1}{36}}\|\psi_\lambda\|_{L^2(\IS^2)},$$
which also yields the result for $2<p\leq 4$ by interpolation. Hence, in order to prove Theorem~\ref{t:maintheo}, we are left with the proof of Proposition~\ref{p:main-prop} which will be the object of the rest of the article.

%As above, we let $0\leq \chi\leq 1$ be a smooth cutoff function which is equal to $1$ on $[-1,1]$ and to $0$ outside $[-2,2]$. Given $x\in \IS^2$, we denote by $\exp_{x}$ the exponential map induced by the metric $g_0$ and we set
%$$\chi_{x,r}(y):=\chi\left(\frac{\|\exp_{x}^{-1}(y)\|}{r}\right)\in\mathcal{C}^{\infty}(\IS^2).$$
%This function is equal to $1$ on $B_r(x)$ and to $0$ outside $B_{2r}(x)$. Hence, 

Coming back to Proposition~\ref{p:main-prop}, it is in fact sufficient to get an uniform upper bound on
$$\tilde{M}_{B_{r_0}(x_0),\alpha,r}(\psi_\lambda):=\sup\left\{\frac{1}{r^{1+\alpha}}\int_{\IS^2}\chi_{x,r}(y)|\psi_\lambda(y)|^2d\upsilon_{g_0}(y):x\in B_{r_0}(x_0)\right\},$$
where we used the conventions of \S\ref{s:kakeya} for the function $\chi_{x,r}$. In order to get this uniform control, we will make use of the invariance properties of semiclassical Wigner distributions that we recently obtained with Maci\`a~\cite{MR16, MR19}. In order to make use of semiclassical methods~\cite{Zw12}, we set $h=\lambda^{-1}$ and $u_h=\psi_{\lambda}$. Hence, one has
\begin{equation}\label{e:semiclassical-schrodinger}
 -h^2\Delta_{g_0}u_h+h^2 Vu_h=u_h,\quad\|u_h\|_{L^2(\IS^2)}=1.
\end{equation}
Let now $x$ be a point in $B_{r_0}(x_0)$ and $h^{\beta}\leq r\leq \frac{\pi}{4}$. In terms of pseudodifferential operators on $\IS^2$~\cite[\S14.2]{Zw12}, the quantity we are interested in can be rewritten as 
$$\int_{\IS^2}\chi_{x,r}(y)|u_h(y)|^2d\upsilon_{g_0}(y)=\left\langle \Op_h\left(\chi_{x,r}\right)u_h,u_h\right\rangle_{L^2(\mathbb{S}^2)},$$
where $\Op_h$ is a semiclassical quantization~\cite[\S14.2.3]{Zw12}. Note that, in order to have $\chi_{x,r}$ amenable to semiclassical pseudodifferential calculus~\cite[\S4.4.1]{Zw12} (see also~\cite[\S2.2, App.A]{DJN19} for the case of manifolds), we need to impose that
\begin{equation}\label{e:admissible}r\geq h^{\beta}\quad\text{and}\quad 0\leq\beta<\frac{1}{2}.\end{equation}
We will now revisit the arguments of~\cite{MR16, MR19} in that specific framework and show how they yield the expected result.

\subsection{Spectral cutoff}\label{ss:cutoff} We fix some smooth cutoff function $0\leq\chi_0\leq 1$ which is equal to $1$ on the interval $[1/2,2]$ and to $0$ outside $[1/4,4]$. Thanks to~\eqref{e:semiclassical-schrodinger}, one has
$$\left\langle \Op_h\left(\chi_{x,r}\right)u_h,u_h\right\rangle_{L^2(\mathbb{S}^2)}=\left\langle \Op_h\left(\chi_{x,r}\right)\chi_0(-h^2\Delta_{g_0}+h^2V)u_h,u_h\right\rangle_{L^2(\mathbb{S}^2)}.$$
According to~\cite[Th.~14.9]{Zw12}, $\chi_0(-h^2\Delta_{g_0}+h^2V)$ is a semiclassical pseudodifferential operator in the class $\Psi^{-\infty}(\IS^2)$ with principal symbol equal to $\chi_0(\|\eta\|_{g_0^*(y)}^2).$ Hence, the composition rule for pseudodifferential operators~\cite[Th.~4.18 and~14.1]{Zw12} implies that
$$\left\langle \Op_h\left(\chi_{x,r}\right)u_h,u_h\right\rangle_{L^2(\mathbb{S}^2)}=\left\langle \Op_h\left(\chi_{x,r}(y)\chi_0(\|\eta\|^2)\right)u_h,u_h\right\rangle_{L^2(\mathbb{S}^2)}+\mathcal{O}(h^{1-2\beta}),$$
where the constant in the remainder is uniform for $x\in \IS^2$ and $r\geq h^\beta$. In the following, we set
$$a_{x,r}(y,\eta):=\chi_{x,r}(y)\chi_0(\|\eta\|_{g_0^*(y)}^2).$$

\subsection{Applying the evolution by the free Schr\"odinger flow}\label{ss:egorov} We write
\begin{equation}\label{e:decompose-Laplace}
 -\Delta_{g_0}= A^2-\frac{1}{4},
\end{equation}
where $A$ is a selfadjoint pseudodifferential operator of order $1$ with principal symbol $\|\eta\|_{g_0^*(y)}$ and satisfying
\begin{equation}\label{e:period-quantum}
e^{2i\pi A}=-\text{Id}.
\end{equation}
Equivalently, one has $A=\sqrt{\frac{1}{4}-\Delta_{g_0}}.$ The eigenvalue equation~\eqref{e:semiclassical-schrodinger} can be rewritten as
$$\left(A^2-\frac{1}{h^2}\right)u_h=\left(\frac{1}{4}-V\right)u_h\quad\Longrightarrow\quad \left(A-\frac{1}{h}\right)u_h=\mathcal{O}_{L^2}(h).$$
In particular, one has
\begin{equation}\label{e:small-quasimode}e^{is \left(A-\frac{1}{h}\right)}u_h=u_h+\int_0^se^{i\tau \left(A-\frac{1}{h}\right)}\left(A-\frac{1}{h}\right)u_hd\tau=u_h+\mathcal{O}_{L^2}(|s|h).\end{equation}
This leads to
\begin{equation}\label{e:av-free-schr}
 \int_{\IS^2}\chi_{x,r}(y)|u_h(y)|^2d\upsilon_{g_0}(y)=\left\langle \left(\frac{1}{2\pi}\int_0^{2\pi}e^{isA}\Op_h\left(a_{x,r}\right)e^{-isA} ds\right)u_h,u_h\right\rangle_{L^2(\mathbb{S}^2)}+\mathcal{O}(h^{1-2\beta})
\end{equation}
In the following, given $a$ in $\ml{C}^{\infty}_c(T^*\IS^2\setminus \underline{0})$, we set, by analogy with the Radon transfom,
$$\ml{R}_{\text{qu}}(\Op_h(a)):=\frac{1}{2\pi}\int_0^{2\pi}e^{is A}\Op_h(a)e^{-isA}ds.$$
According to Remark~\ref{r:egorov} below, the Egorov Theorem allows to relate the operator $\ml{R}_{\text{qu}}(\Op_h(a_{x,r}))$ to the classical average by the geodesic flow:
\begin{equation}\label{e:egorov}
 \ml{R}_{\text{qu}}(\Op_h(a_{x,r}))= \Op_h\left(\frac{1}{2\pi}\int_0^{2\pi} a_{x,r}\circ \varphi_0^{t}dt\right)+\mathcal{O}_{L^2\rightarrow L^2}(h^{1-2\beta}),
\end{equation}
where the constant in the remainder is uniform for $x\in\IS^2$ and $r\geq h^{\beta}$ and where $\varphi_0^t$ is the Hamiltonian flow associated with the Hamiltonian function\footnote{This is just a reparametrization of the standard geodesic flow.} $H_0(y,\eta):=\|\eta\|_{g_0(y)}.$ Given $a$ in $\ml{C}^{\infty}_c(T^*\IS^2\setminus \underline{0})$, we set 
$$\ml{R}_{\text{cl}}(a):=\frac{1}{2\pi}\int_0^{2\pi} a\circ \varphi_0^{t}dt.$$

\begin{rema}\label{r:egorov} Let us briefly remind how to prove~\eqref{e:egorov}. This is standard~\cite[App.~A.3]{DJN19} and we just need to pay attention to our class of symbols. First, we write, for every $s,t\in[0,2\pi]$,
 $$\frac{d}{ds}\left(e^{is A}\Op_h(a_{x,r}\circ\varphi_0^{t-s})e^{-isA}\right)=e^{is A}\left(\frac{i}{h}\left[hA,\Op_h(a_{x,r}\circ\varphi_0^{t-s})\right]-\Op_h\left(\{H_0,a_{x,r}\circ\varphi_0^{t-s}\}\right)\right)e^{-isA}.$$
 We now let $\chi_1$ be a smooth function which is equal to $1$ in a neighborhood of $[1/4,4]$ and to $0$ outside $[1/8,8]$. In particular, $\chi_1(H_0^2)$ is equal to $1$ on the support of $a_{x,r}$. Combining this with the composition rules for pseudodifferential operators with exotic symbols on manifolds~\cite[Lemma~A.6]{DJN19}, we know that, for every $\tau\in [0,2\pi]$,
 \begin{eqnarray*}\Op_h(a_{x,r}\circ\varphi_0^\tau)&=&\Op_h(a_{x,r}\circ\varphi_0^\tau)\Op_{h}(\chi_1(H_0^2))+\mathcal{O}_{L^2\rightarrow L^2}(h^2)\\
&=&\Op_{h}(\chi_1(H_0^2))\Op_h(a_{x,r}\circ\varphi_0^{\tau})+\mathcal{O}_{L^2\rightarrow L^2}(h^2).  
 \end{eqnarray*}
 We can also remark using the composition rules for pseudodifferential operators that
 $$hA\Op_{h}(\chi_1(H_0^2))=\Op_h(\chi_1(H_0^2))hA+h\Op_h(r)+\mathcal{O}_{L^2\rightarrow L^2}(h^2),$$
 where $r$ is a smooth compactly supported function that depends in a multilinear way of the derivatives of order $\geq 1$ of the function $\chi_1(H_0^2)$. Thus its support does not intersect the support of $a_{x,r}$. In particular, using the composition rule~\cite[Lemma~A.6]{DJN19} one more time and the support properties of $a_{x,r}$, one has $\Op_h(a_{x,r})\Op_h(r)=\mathcal{O}_{L^2\rightarrow L^2}(h^2)$.
Hence, after integration over the interval $[0,2\pi]$ and applying the Calder\'on-Vaillancourt Theorem, one finds
 \begin{eqnarray*}
   \ml{R}_{\text{qu}}(\Op_h(a_{x,r})) &= &\Op_h\left(\frac{1}{2\pi}\int_0^{2\pi} a_{x,r}\circ \varphi_0^{t}dt\right) +\mathcal{O}_{L^2\rightarrow L^2}(h)\\  
  &+&\frac{1}{2\pi }\int_0^{2\pi}\int_0^t\left(\frac{i}{h}\left[hA\Op_{h}(\chi_1(H_0^2)),\Op_h(a_{x,r}\circ\varphi_0^{t-s})\right]\right) dsdt\\
  &-&\frac{1}{2\pi }\int_0^{2\pi}\int_0^t\Op_h\left(\{H_0,a_{x,r}\circ\varphi_0^{t-s}\}\right) dsdt.
 \end{eqnarray*}
%where $\tilde{u}_h(s)=e^{-isA} \tilde{u}_h$. %We now let $\chi_1$ be a smooth function which is equal to $1$ in a neighborhood of $[1/4,4]$ and to $0$ outside $[1/8,8]$. In particular, $\chi_1(H_0^2)$ is equal to $1$ on the support of $a_{x,r}$ and 
%$$u_h=\chi_1(-h^2\Delta_g+h^2 V)u_h=\chi_1(-h^2\Delta_g)u_h +\mathcal{O}_{L^2}(h^2).$$
%In particular, using the Calder\'on-Vaillancourt Theorem for pseudodifferential operators, we can rewrite the above equality as
% \begin{eqnarray*}
%  \left\langle \ml{R}_{\text{qu}}(\Op_h(a_{x,r})) u_h,u_h\right\rangle &= &\left\langle \Op_h\left(\ml{R}_{\text{cl}}(a_{x,r})\right)u_h,u_h\right\rangle +\mathcal{O}(h)\\  
%  &+&\frac{1}{2\pi }\int_0^{2\pi}\int_0^t\left(\frac{i}{h}\left\langle\left[hA\chi_1^2(-h^2\Delta_{g_0}),\Op_h(a_{x,r}\circ\varphi_0^{t-s})\right]u_h(s),u_h(s)\right\rangle\right) dsdt\\
%  &-&\frac{1}{2\pi }\int_0^{2\pi}\int_0^t\left\langle\Op_h\left(\{H_0,a_{x,r}\circ\varphi_0^{t-s}\}\right)u_h(s),u_h(s)\right\rangle dsdt,
% \end{eqnarray*}
% where $u_h(s)=e^{-isA} u_h$ and where the constant in the remainder is uniform for $x\in \mathbb{S}^2$ and $r\geq h^{\beta}$. 
As all our pseudodifferential operators are microlocally supported in a compact\footnote{This was the main reason for inserting the pseudodifferential cutoff $\Op_{h}(\chi_1(H_0^2))$.} set of $T^*\IS^2$, we can again apply the composition rule for exotic symbols on a compact manifold as stated in~\cite[Lemma~A.6]{DJN19}. Thus, we can conclude that~\eqref{e:egorov} holds. Inspecting carefully the argument, we can in fact conclude that
\begin{lemm} With the above conventions, one can find $\tilde{a}_{x,r}\in S^{\text{comp}}_{\beta}(T^*\IS^2)$ (as defined in~\cite[\S 2.2]{DJN19}) such that
\begin{equation}\label{e:egorov2}\ml{R}_{\text{qu}}(\Op_h(a_{x,r}))=\Op_h(\tilde{a}_{x,r})+\mathcal{O}_{L^2\rightarrow L^2} (h^2),\end{equation}
 where the constant in the remainder is uniform for $x\in\IS^2$ and $r\geq h^{\beta}$. Moreover, $\tilde{a}_{x,r}$ is equal to $\mathcal{R}_{\text{cl}}(a_{x,r})$ modulo $h^{1-2\beta}S^{\text{comp}}_{\beta}(T^*\IS^2)$ and its support is contained in the support of $\mathcal{R}_{\text{cl}}(a_{x,r})$.
\end{lemm}

%where the constant in the remainder is also uniform uniform for $x\in\IS^2$ and $r\geq h^{\beta}$.

%. The same argument would in fact shows that
%  \begin{equation}\label{e:egorov2}\chi_1(-h^2\Delta_{g_0})\left(\ml{R}_{\text{qu}}(\Op_h(a_{x,r}))-\Op_h\left(\ml{R}_{\text{cl}}(a_{x,r})\right)\right)\chi_1(-h^2\Delta_{g_0})=h^{1-2\beta}\Op_h\left(R_{x,r}\right)+\mathcal{O}_{L^2\rightarrow L^2}(h^2),\end{equation}
%  where $R_{x,r}$ is a symbol in the exotic class of symbols $S^{\text{comp}}_{\beta}(T^*\mathbb{S}^2)$~\cite[\S 2.2.1]{DJN19} (with all seminorms that can be made uniform for $x\in\IS^2$ and $r\geq h^{\beta}$) and where the constant in the remainder is also uniform in these parameters.
\end{rema}

\begin{rema}\label{r:semiclassical-pert} The arguments used from the beginning of this Section would work as well for the following semiclassical problem:
$$-h^2\Delta_{g_0}u_h+\varepsilon_h Vu_h=u_h,\quad\|u_h\|_{L^2(\IS^2)}=1,$$
 where $\varepsilon_h\rightarrow 0$ fast enough. More precisely, the above proofs only require $h^{-1}\varepsilon_h\rightarrow 0$ in order to have a small remainder in~\eqref{e:small-quasimode}. In this case, this would yield the bound
 $$\int_{\IS^2}\chi_{x,r}(y)|u_h(y)|^2d\upsilon_{g_0}(y)=\left\langle \left(\frac{1}{2\pi}\int_0^{2\pi}e^{-isA}\Op_h\left(a_{x,r}\right)e^{isA} ds\right)u_h,u_h\right\rangle_{L^2(\mathbb{S}^2)}+\mathcal{O}(h^{1-2\beta})+\ml{O}(h^{-1}\varepsilon_h).$$
 The argument from~\cite{MR16} would allow to remove this extra remainder $\ml{O}(h^{-1}\varepsilon_h)$ and to handle the case $\varepsilon_h\rightarrow 0^+$. Yet, as this kind of condition on the size of the potential already appeared in Remarks~\ref{r:Lp-semiclassical} and~\ref{r:KN-Lp}, we do not pursue this here.
\end{rema}

\subsection{Weinstein averaging method} 

Following Weinstein~\cite{W77}, one can use~\eqref{e:period-quantum} to obtain the following exact commutation relation:
$$\left[\ml{R}_{\text{qu}}(\Op_h(a_{x,r})),A\right]=0.$$
In particular, thanks to~\eqref{e:decompose-Laplace}, one has
\begin{equation}\label{e:commute}
\left[\ml{R}_{\text{qu}}(\Op_h(a_{x,r})),\Delta_{g_0}\right]=0.
\end{equation}
Using~\eqref{e:semiclassical-schrodinger}, this implies that
$$\left\langle \left[V,\ml{R}_{\text{qu}}(\Op_h(a_{x,r}))\right]u_h,u_h\right\rangle_{L^2(\mathbb{S}^2)}=0.$$
Thanks to~\eqref{e:egorov2}, this can be rewritten as
$$\left\langle \left[V,\Op_h(\tilde{a}_{x,r})\right]u_h,u_h\right\rangle_{L^2(\mathbb{S}^2)}=\mathcal{O}(h^2).$$
%As in Remark~\ref{r:egorov}, we set $\tilde{u}_h=\chi_1(-h^2\Delta_g)u_h.$ Using the composition rules and the functional calculus for pseudodifferential operators, one has 
%$$(-h^2\Delta_{g_0}+h^2V)\chi_1(-h^2\Delta_{g_0})=\chi_1(-h^2\Delta_{g_0})(-h^2\Delta_{g_0}+h^2V)+\mathcal{O}_{L^2\rightarrow L^2}(h^3).$$
%In particular, $(-h^2\Delta_{g_0}+h^2V)\tilde{u}_h=\tilde{u}_h+\mathcal{O}(h^3)$ and
%$$\left\langle \left[V,\ml{R}_{\text{qu}}(\Op_h(a_{x,r}))\right]u_h,u_h\right\rangle_{L^2(\mathbb{S}^2)}=0.$$
As in Remark~\ref{r:egorov}, we can insert pseudodifferential cutoffs and we find
$$\left\langle \left[V\Op_h(\chi_1(H_0^2)),\Op_h(\tilde{a}_{x,r})\right]u_h,u_h\right\rangle_{L^2(\mathbb{S}^2)}=\mathcal{O}(h^2).$$
Hence, thanks to the composition rule for pseudodifferential operators~\cite[Lemma~A.6]{DJN19} with exotic symbols, we get 
$$\left\langle \Op_h\left(\left\{V,\mathcal{R}_{\text{cl}}(a_{x,r})\right\}\right)u_h,u_h\right\rangle_{L^2(\mathbb{S}^2)}=\ml{O}(h^{1-3\beta}),$$
where the constant in the remainder is uniform for $x\in \IS^2$ and $r\geq h^\beta$. Observe that the extra loss in $\ml{O}(h^{1-3\beta})$ (compared with $\ml{O}(h^{1-2\beta})$) comes from the subprincipal term in $\tilde{a}_{x,r}$. Applying the argument of paragraph~\ref{ss:egorov} one more time, we find that
$$\left\langle \Op_h\left(\frac{1}{2\pi}\int_0^{2\pi}\left\{V,\mathcal{R}_{\text{cl}}(a_{x,r})\right\}\circ\varphi_0^tdt\right)u_h,u_h\right\rangle_{L^2(\mathbb{S}^2)}=\ml{O}(h^{1-3\beta}),$$
from which we infer
$$\left\langle \Op_h\left(\left\{\ml{R}_{\text{cl}}(V),\mathcal{R}_{\text{cl}}(a_{x,r})\right\}\right)u_h,u_h\right\rangle_{L^2(\mathbb{S}^2)}=\ml{O}(h^{1-3\beta}),$$
with the constant in the remainder enjoying the same uniformity property as before. Here $V$ is identified with its pullback on $T^*\IS^2\setminus \underline{0}$ via the canonical projection $\Pi(y,\eta)=y$. 

Let us now denote by $\varphi^t_{\langle V\rangle}$ the Hamiltonian flow induced by $\ml{R}_{\text{cl}}(V)$. As $\ml{R}_{\text{cl}}(V)$ and $H_0$ Poisson commute, one has $\varphi_0^t\circ\varphi_{\langle V\rangle}^s=\varphi^s_{\langle V\rangle}\circ\varphi_0^t$ for every $t$ and $s$ in $\mathbb{R}$. We note that all the above argument would work as well if we replace $a_{x,r}$ by $a_{x,r}\circ\varphi_{\langle V\rangle}^\tau$ and the remainder would remain uniform in $\tau$ (and in $(x,r)$) provided that $\tau$ remains on a bounded interval. Hence, one has, uniformly for $\tau\in[-\tau_0,\tau_0]$, $x\in \IS^2$ and $r\geq h^{\beta}$,
\begin{equation}\label{e:inv-V}\left\langle \Op_h\left(\left\{\ml{R}_{\text{cl}}(V),\mathcal{R}_{\text{cl}}(a_{x,r})\circ\varphi_{\langle V\rangle}^{\tau}\right\}\right)u_h,u_h\right\rangle_{L^2(\mathbb{S}^2)}=\ml{O}(h^{1-3\beta}).
\end{equation}
We integrate this expression between $0$ and $\tau$:
$$\left\langle \Op_h\left(\mathcal{R}_{\text{cl}}(a_{x,r})\circ\varphi_{\langle V\rangle}^{\tau}\right)u_h,u_h\right\rangle_{L^2(\mathbb{S}^2)}=\left\langle \Op_h\left(\mathcal{R}_{\text{cl}}(a_{x,r})\right)u_h,u_h\right\rangle_{L^2(\mathbb{S}^2)}+\ml{O}(h^{1-3\beta}).$$
Combining this with~\eqref{e:av-free-schr}, we find
\begin{equation}\label{e:final-average}\int_{\IS^2}\chi_{x,r}(y)|u_h(y)|^2d\upsilon_{g_0}(y)=\left\langle \Op_h\left(\frac{1}{2\tau_0}\int_{-\tau_0}^{\tau_0}\mathcal{R}_{\text{cl}}(a_{x,r})\circ\varphi_{\langle V\rangle}^{\tau}d\tau\right)u_h,u_h\right\rangle_{L^2(\mathbb{S}^2)}+\mathcal{O}(h^{1-3\beta}),\end{equation}
where the constant in the remainder is uniform for $x$ in $K$ and $r\geq h^{\beta}$.
\begin{rema} Rather than for studying eigenfunctions, Weinstein's argument was initially developed to study the distribution of eigenvalues of $-\Delta_{g_0}+V$ inside each cluster near $\lambda_j^2=j(j+1)$~\cite{W77}. This was achieved by showing via this kind of averaging arguments that $-\Delta_{g_0}+V$ is conjugated to $-\Delta_{g_0}+\ml{R}_{\text{qu}}(V)$ modulo small error terms. See~\cite{CdV79, Gu78, Gu81, Ze96, Ze97} for further developments on these eigenvalue problems.           
            \end{rema}

\subsection{Applying Calder\'on-Vaillancourt Theorem}

We are now in position to apply the Calder\'on-Vaillancourt Theorem~\cite[Th.~5.1]{Zw12} which tells us that
$$\left\|\Op_h\left(\frac{1}{2\tau_0}\int_{-\tau_0}^{\tau_0}\mathcal{R}_{\text{cl}}(a_{x,r})\circ\varphi_{\langle V\rangle}^{\tau}d\tau\right)\right\|_{L^2\rightarrow L^2}\leq C\left\|\frac{1}{2\tau_0}\int_{-\tau_0}^{\tau_0}\mathcal{R}_{\text{cl}}(a_{x,r})\circ\varphi_{\langle V\rangle}^{\tau}d\tau\right\|_{L^{\infty}(T^*\mathbb{S}^2)}+\mathcal{O}(h^{1-3\beta}),$$
where $C_0$ is some universal constant and where the constant in the remainder is one more time uniform for $x$ in $\IS^2$ and $r\geq h^{\beta}$. Together with~\eqref{e:final-average}, we finally get
$$\int_{\IS^2}\chi_{x,r}(y)|u_h(y)|^2d\upsilon_{g_0}(y)\leq C\left\|\frac{1}{2\tau_0}\int_{-\tau_0}^{\tau_0}\mathcal{R}_{\text{cl}}(a_{x,r})\circ\varphi_{\langle V\rangle}^{\tau}d\tau\right\|_{L^{\infty}(T^*\mathbb{S}^2)}+\mathcal{O}(h^{1-3\beta}).$$ 
From the construction of $a_{x,r}$, one can in fact reduce to the unit cotangent bundle and conclude that the following key lemma holds
\begin{lemm} With the above conventions, one has
\begin{equation}\label{e:final-bound-ball}\int_{\IS^2}\chi_{x,r}(y)|u_h(y)|^2d\upsilon_{g_0}(y)\leq C\left\|\frac{1}{4\pi\tau_0}\int_{-\tau_0}^{\tau_0}\int_0^{2\pi}\chi_{x,r}\circ\varphi^t_0\circ\varphi_{\langle V\rangle}^{\tau}dtd\tau\right\|_{L^{\infty}(S^*\mathbb{S}^2)}+\mathcal{O}(h^{1-3\beta}),\end{equation}
where we identify $\chi_{x,r}$ with its pullback on $S^*\IS^2$ and where the constant in the remainder is uniform for $x$ in $\IS^2$ and $r\geq h^{\beta}$.
\end{lemm}

In order to facilitate the discussion, we shall work on the space of geodesic $G(\mathbb{S}^2)\simeq\IS^2$. With the induced symplectic form on $\IS^2$, $\varphi^{\tau}_{\langle V\rangle}$ can be viewed as the Hamiltonian flow of $\ml{R}(V)$ on $\IS^2$. Hence, what we are aiming at is an upper bound on
$$0\leq \frac{1}{2\tau_0}\int_{-\tau_0}^{\tau_0}\mathcal{R}(\chi_{x,r})\circ\varphi_{\langle V\rangle}^{\tau}(\gamma)d\tau,$$
when $\gamma\in G(\IS^2)\simeq\IS^2$ and when $r\ll\tau_0$. It is in fact sufficient to find an upper bound on
$$\frac{1}{2\tau_0}\int_{-\tau_0}^{\tau_0}\mathcal{R}(\mathbf{1}_{B_{2r}(x)})\circ\varphi_{\langle V\rangle}^{\tau}(\gamma)d\tau,$$
where $\mathbf{1}_{B_{2r}(x)}$ is the characteristic function of the geodesic ball of radius $2r$ centered at $x$. The function $\mathcal{R}(\mathbf{1}_{B_{2r}(x)})$ is supported in a neighborhood of width $4r$ of $\Gamma_x\subset G(\IS^2)$ and it is bounded from above by $4r$. Hence,
\begin{equation}\label{e:trivialbound}\forall\gamma\in G(\IS^2),\quad0\leq \frac{1}{2\tau_0}\int_{-\tau_0}^{\tau_0}\mathcal{R}(\mathbf{1}_{B_{2r}(x)})\circ\varphi_{\langle V\rangle}^{\tau}(\gamma)d\tau\leq 4r.\end{equation}

\begin{rema}\label{r:semiclassical-pert-2} In the case of semiclassical Schr\"odinger operators as in Remark~\ref{r:semiclassical-pert}, the argument would work similarly and we would also obtain the bound~\eqref{e:final-bound-ball} for this semiclassical problem (up to the already extra remainder $\mathcal{O}(h^{-1}\varepsilon_h)$ that apeared in this Remark).
\end{rema}

\subsection{Flow lines of $\varphi_{\langle V\rangle}^t$ near $\Gamma_{x_0}$}

So far we did not use our assumptions on $V$ or on the point $x_0$. They will now be used to get an improvement of order $r^{1/2}$ on the upper bound~\eqref{e:trivialbound} when $x\in B_{r_0}(x_0)$. To that aim, we now fix $x_0$ satisfying the assumption of the Theorem and we will analyze the flow lines of $\varphi_{\langle V\rangle}^t$ near a given point $\gamma_0$ of $\Gamma_{x_0}$. 

Without loss of generality, we may suppose that $x_0$ is the north pole, i.e. with coordinates $(0,0,1)$ in the representation~\eqref{e:sphere}. Then, for every $x\in B_{\epsilon_0}(x_0)$, $\Gamma_x$ is a great circle of the sphere lying in the annulus
$$\mathcal{A}_{\epsilon_0}:=\left\{(x_1,x_2,x_3)\in\mathbb{R}^3:x_1^2+x_2^2+x_3^2=1,\ |x_3|\leq\sin\epsilon_0\right\}.$$
Similarly, the function $\mathcal{R}(\mathbf{1}_{B_{2r}(x)})$ is supported on an annulus of width $2|\sin(2r)|$ around $\Gamma_x$ and it takes the value $4r$ on this annulus. In particular, if $\tau_0>0$ and $r_1>0$ are chosen small enough, then, for every $x\in B_{\epsilon_0}(x_0)$ and for every $0<r<r_1$, the support of 
\begin{equation}\label{e:av-radon-potential}\frac{1}{2\tau_0}\int_{-\tau_0}^{\tau_0}\mathcal{R}(\mathbf{1}_{B_{2r}(x)})\circ\varphi_{\langle V\rangle}^{\tau}d\tau\end{equation}
is contained in the annulus $\mathcal{A}_{2\epsilon_0}$. Hence, once we have fixed $x\in B_{\epsilon_0}(x_0)$, we just need to study the value of this function inside such an annulus. More precisely, we want to show that this is of order $\mathcal{O}(r^{3/2})$ uniformly for $\gamma$ in this annulus.

%All these annuli are naturally oriented thanks to the orientation on $\IS^2$ and this induces an orientation of their (smooth) boundaries. %We also endow each $\Gamma_x$ with the same orientation as these boundaries curves.

Let $\gamma_0\in \Gamma_{x_0}$ and let us prove this uppper bound in a neighborhood of a fixed $\gamma_0$. Without loss of generality, we can suppose that, in spherical coordinates $(\phi,\theta)$, one has $\gamma_0=(\pi/2,0)$. The vector field $X_{\langle V\rangle}$ can be written in this system of coordinates:
$$X_{\langle V\rangle}(\phi,\theta)=-\frac{1}{\sin\phi}\frac{\partial\mathcal{R}(V)}{\partial\theta}\partial_\phi+\frac{\partial\mathcal{R}(V)}{\partial\phi}\partial_\theta.$$
We need to distinguish two situations: 
\begin{enumerate}
 \item $X_{\langle V\rangle}(\gamma_0)\notin T_{\gamma_0}\Gamma_{x_0}$ which means that $\frac{\partial\mathcal{R}(V)}{\partial\theta}(\pi/2,0)\neq 0$;
 \item $X_{\langle V\rangle}(\gamma_0)\in T_{\gamma_0}\Gamma_{x_0}$ which means that $\frac{\partial\mathcal{R}(V)}{\partial\theta}(\pi/2,0)= 0$. In that case, the hypothesis of Theorem~\ref{t:maintheo} implies that $\frac{\partial\mathcal{R}(V)}{\partial\phi}(\pi/2,0)\neq 0$ and $\frac{\partial^2\mathcal{R}(V)}{\partial\theta^2}(\pi/2,0)\neq 0$
\end{enumerate}
The Hamilton-Jacobi equations can be written as
\begin{equation}\label{e:HJ}
 \phi'(\tau)=-\frac{1}{\sin\phi(\tau)}\frac{\partial\mathcal{R}(V)}{\partial\theta}(\phi(\tau),\theta(\tau)),\quad\text{and}\quad\theta'(\tau)=\frac{\partial\mathcal{R}(V)}{\partial\phi}(\phi(\tau),\theta(\tau)).
\end{equation}

\subsubsection{The transverse case}

Let us begin with the first situation which is slightly easier to handle. Witout loss of generality, we can suppose that $\frac{\partial\mathcal{R}(V)}{\partial\theta}(\pi/2,0)>0$ (the negative case is handled similarly). First, using spherical coordinates, we fix an open neighborhood $\mathcal{U}_{2\epsilon_0}:=(\pi/2-4\epsilon_0,\pi/2+4\epsilon_0)\times(-2\epsilon_0,2\epsilon_0)$ so that
\begin{equation}\label{e:transverse-local}
\forall \gamma=(\phi,\theta)\in\mathcal{U}_{2\epsilon_0},\quad \frac{\partial\mathcal{R}(V)}{\partial\theta}(\phi,\theta)>\frac{1}{2}\frac{\partial\mathcal{R}(V)}{\partial\theta}(\pi/2,0)=:a_0>0.
\end{equation}
Up to decreasing the value $\tau_0$, we can suppose without loss of generality that $\varphi^\tau_{\langle V\rangle}(\gamma)$ belongs to $\mathcal{U}_{2\epsilon_0}$ for every $|\tau|\leq\tau_0$ and for every $\gamma\in\mathcal{U}_{\epsilon_0}$. As already explained, the support of~\eqref{e:av-radon-potential} is contained in $\mathcal{A}_{2\epsilon_0}$. For the moment, we will study locally its value inside $\mathcal{U}_{\epsilon_0}\subset\mathcal{A}_{2\epsilon_0}$. We now fix some $\gamma$ in $\mathcal{U}_{\epsilon_0}$. In particular,
$$\forall |\tau|\leq \tau_0,\quad  \frac{\partial\mathcal{R}(V)}{\partial\theta}\left(\varphi^\tau_{\langle V\rangle}(\gamma)\right)\geq a_0,$$
which implies thanks to~\eqref{e:HJ} that $\phi'(\tau)<0$ along this piece of trajectory. This yields the following upper bound along the orbit $\left(\varphi^\tau_{\langle V\rangle}(\gamma)\right)_{-\tau_0\leq\tau\leq\tau_0}$:
\begin{equation}\label{e:angle-interval}\phi(\tau_2)-\phi(\tau_1)\leq -\frac{a_0}{\cos(4\epsilon_0)}(\tau_2-\tau_1)\ \Longleftrightarrow\  \tau_2-\tau_1\leq\frac{\cos(4\epsilon_0)}{a_0}(\phi(\tau_1)-\phi(\tau_2)),\end{equation}
for every $-\tau_0\leq\tau_1\leq\tau_2\leq\tau_0$. 

Recall now that the function in~\eqref{e:av-radon-potential} is defined by averaging $\mathcal{R}(\mathbf{1}_{B_{2r}(x)})$ for some $x\in B_{\epsilon_0}(x_0)$ and some $0<r<r_1$. In spherical coordinates, $x$ can be written $(\phi_x,\theta_x)$ where $0\leq\phi_x\leq \epsilon_0$ and $0\leq\theta_x\leq 2\pi$. Hence, using our identification $G(\IS^2)\simeq\IS^2$, $\mathcal{R}(\mathbf{1}_{B_{2r}(x)})$ is $4r$ times the characteristic function of the annulus of width $4r$ centered at $\Gamma_x$,
$$\mathcal{A}_{2r}(x)=\left\{(\phi,\theta):\ \phi-\text{arccos}\left(-\cos(\theta-\theta_x)\sin(\phi_x)\right)\in[-2r,2r],\ 0\leq \theta\leq 2\pi\right\}.$$
The boundary of this annulus is given by
$$\partial \mathcal{A}_{2r}(x)=\left\{\left(\text{arccos}\left(-\cos(\theta-\theta_x)\sin(\phi_x)\right)\pm 2r,\theta\right):\ 0\leq\theta\leq 2\pi\right\}$$
and it is oriented thanks to the natural orientation on $\IS^2$. Using now that $\mathcal{R}(V)$ is of class $\mathcal{C}^1$ and~\eqref{e:transverse-local}, we know that, up to decreasing the value of $\epsilon_0$ (and thus of $\tau_0$ and $r_1$), the vector field $X_{\langle V\rangle}$ is uniformly (negatively) transverse to $\partial \mathcal{A}_{2r}(x)\cap\mathcal{U}_{2\epsilon_0}$ for every $x\in B_{\epsilon_0}(x_0)$ and for every $0<r<r_1$. In particular, given $\gamma\in\mathcal{U}_{\epsilon_0}$, the set
$$\left\{\tau\in[-\tau_0,\tau_0]:\varphi^{\tau}_{\langle V\rangle}(\gamma)\in \mathcal{A}_{2r}(x)\right\}$$
is an interval that we denote by $I_{x,r}(\gamma)$. Hence,
$$0\leq\frac{1}{2\tau_0}\int_{-\tau_0}^{\tau_0}\mathcal{R}(\mathbf{1}_{B_{2r}(x)})\circ\varphi_{\langle V\rangle}^{\tau}(\gamma)d\tau\leq \frac{2r|I_{x,r}(\gamma)|}{\tau_0},$$
and it remains to determine an upper bound on the size of this interval in terms of $r$. Thanks to the upper bound~\eqref{e:angle-interval}, the length of the interval is bounded by the maximal variation of $\phi$ along the orbit of $\gamma$ inside $\mathcal{A}_{2r}(x)$. If we denote the interval $I_{x,r}(\gamma)$ by $[\tau_1,\tau_2]$, then
$$\phi(\tau_1)-\phi(\tau_2)\leq 4r+\left|\text{arccos}\left(-\cos(\theta(\tau_1)-\theta_x)\sin(\phi_x)\right)-\text{arccos}\left(-\cos(\theta(\tau_2)-\theta_x)\sin(\phi_x)\right)\right|.$$
As $\phi_x\in[-\epsilon_0,\epsilon_0]$ (with $\epsilon_0>0$ small), this yields an upper bound of the form
$$\phi(\tau_1)-\phi(\tau_2)\leq 4r+C \sin(\epsilon_0)|\tau_2-\tau_1|,$$
where $C>0$ is some uniform constant. Combined with~\eqref{e:angle-interval}, it gives us
$$0\leq|I_{x,r}(\gamma)|=\tau_2-\tau_1\leq\frac{4r}{1-C\sin(\epsilon_0)},$$
and then, for every $x\in B_{\epsilon_0}(x_0)$ and every $0\leq r\leq r_1$,
\begin{equation}\label{e:bound-transverse-case}\forall\gamma\in\mathcal{U}_{\epsilon_0},\ 0\leq\frac{1}{2\tau_0}\int_{-\tau_0}^{\tau_0}\mathcal{R}(\mathbf{1}_{B_{2r}(x)})\circ\varphi_{\langle V\rangle}^{\tau}(\gamma)d\tau\leq \frac{8r^2}{\tau_0(1-C\sin(\epsilon_0))}.\end{equation}
This shows the expected upper bound in the neighborhood $\mathcal{U}_\epsilon(\gamma_0):=\mathcal{U}_{\epsilon_0}$ of $\gamma_0$ when $X_{\langle V\rangle}(\gamma_0)$ is transverse to $\Gamma_{x_0}$.

\subsubsection{The tangent case}

We now deal with the slightly more delicate case where $X_{\langle V\rangle}(\gamma_0)$ is tangent to $\Gamma_{x_0}$ where $\frac{\partial\mathcal{R}(V)}{\partial\theta}(\pi/2,0)=0$. Thanks to~\eqref{e:hyp-crit}, we can again without loss of generality assume that $\frac{\partial\mathcal{R}(V)}{\partial\phi}(\pi/2,0)>0$, and suppose that 
\begin{equation}\label{e:tangent-local}
\forall \gamma=(\phi,\theta)\in\mathcal{U}_{2\epsilon_0},\quad \frac{\partial\mathcal{R}(V)}{\partial\phi}(\phi,\theta)>\frac{1}{2}\frac{\partial\mathcal{R}(V)}{\partial\phi}(\pi/2,0)=:a_0>0.
\end{equation}
Moreover, thanks to hypothesis~\eqref{e:hyp-trans}, the critical point at $0$ of the map $\theta\mapsto \mathcal{R}(V)(\pi/2,\theta)$ is nondegenerate. In particular, without loss of generality and up to decreasing the value of $\epsilon_0$, there exists $b_0>0$ such that
\begin{equation}\label{e:tangent-local-second}
\forall \gamma=(\phi,\theta)\in\mathcal{U}_{2\epsilon_0},\quad \frac{\partial^2\mathcal{R}(V)}{\partial\theta^2}(\phi,\theta)>\frac{1}{2}\frac{\partial^2\mathcal{R}(V)}{\partial\theta^2}(\pi/2,0)=:b_0>0.
\end{equation}
We now fix $\gamma\in\mathcal{U}_{\epsilon_0}\subset\mathcal{A}_{2\epsilon_0}$ and, as before, we can suppose that, for every $|\tau|\leq \tau_0$,
$$ \frac{\partial\mathcal{R}(V)}{\partial\phi}\left(\varphi^\tau_{\langle V\rangle}(\gamma)\right)\geq a_0\quad \text{and}\quad  \frac{\partial^2\mathcal{R}(V)}{\partial\theta^2}\left(\varphi^\tau_{\langle V\rangle}(\gamma)\right)\geq b_0.$$
As in the transverse case, one has 
$$0\leq\frac{1}{2\tau_0}\int_{-\tau_0}^{\tau_0}\mathcal{R}(\mathbf{1}_{B_{2r}(x)})\circ\varphi_{\langle V\rangle}^{\tau}(\gamma)d\tau\leq \frac{2r|I_{x,r}(\gamma)|}{\tau_0},$$
where
$$I_{x,r}(\gamma):=\left\{\tau\in[-\tau_0,\tau_0]:\varphi^{\tau}_{\langle V\rangle}(\gamma)\in \mathcal{A}_{2r}(x)\right\}.$$

The main difference with the above case is that this set is not an interval in general. Yet, we can note that, along the trajectory of $\gamma$, the vector $(\phi'(\tau),\theta'(\tau))$ is nonvanishing thanks to~\eqref{e:tangent-local}. Moreover, it is tangent to $\partial\mathcal{A}_{r'}(x)$ (for some $r'<2\epsilon_0$) if and only if
$$F(\tau):=\phi'(\tau)-\theta'(t)\frac{\sin(\theta(\tau)-\theta_x)\sin(\phi_x)}{\sqrt{1-\cos^2(\theta(\tau)-\theta_x)\sin^2\phi_x}}=0.$$
We can observe that
\begin{eqnarray*}
F'(\tau)&=&-\frac{1}{\sin\phi(\tau)}\frac{\partial\mathcal{R}(V)}{\partial\theta}(\phi(\tau),\theta(\tau))\frac{\partial^2\mathcal{R}(V)}{\partial\theta\partial\phi}(\phi(\tau),\theta(\tau))\\
&-&\frac{1}{\sin\phi(\tau)}\frac{\partial\mathcal{R}(V)}{\partial\phi}(\phi(\tau),\theta(\tau))\frac{\partial^2\mathcal{R}(V)}{\partial\theta^2}(\phi(\tau),\theta(\tau))+\mathcal{O}(\epsilon_0),
\end{eqnarray*}
where the constant in the remainder is uniformly bounded for $\tau\in[-\tau_0,\tau_0]$ and $\gamma\in\mathcal{U}_{\epsilon_0}$. Thus, as $\partial_\theta\mathcal{R}(V)(\pi/2,0)=0$, we can suppose that, up to decreasing the value of $\epsilon_0>0$, $|F'(\tau)|\geq a_0b_0/2$. In particular, $F$ is monotone and it vanishes at most at one point inside $[-\tau_0,\tau_0]$. As a consequence, the set $I_{x,r}(\gamma)$ is the union of at most two disjoint intervals inside $[-\tau_0,\tau_0]$ that we denote by $[\tau_1,\tau_2]$ and $[\tau_3,\tau_4]$. Moreover, $X_{\langle V\rangle}(\varphi_{\langle V\rangle}^\tau(\gamma))$ is tangent to $\partial A_{2r}(x)$ at most at one point inside $[\tau_1,\tau_2]\cup[\tau_3,\tau_4]$.

It now remains to bound the length of these two intervals in terms of $r$. To that aim, we observe that, for $\tau\in[\tau_1,\tau_2]\cup[\tau_3,\tau_4]$, one can find $r(\tau)\in[-2r,2r]$ such that
$$\phi(\tau)=r(\tau)+\text{arccos}\left(-\cos(\theta(\tau)-\theta_x)\sin\phi_x\right).$$
Given now $\tau,\tau'\in[\tau_1,\tau_2]\cup[\tau_3,\tau_4]$, one finds
$$r(\tau)-r(\tau')=F(\tau')(\tau-\tau')+\frac{F'(\tau')}{2}(\tau-\tau')^2+\mathcal{O}((\tau-\tau')^3),$$
where the constant in the remainder can be made uniform in terms of $r$, $\gamma$, $\tau$ and $\tau'$. We use this equality to find an upper bound on the length of $[\tau_1,\tau_2]$. The other interval (if non empty) is handled simlarly. Recall from the above calculation that $F'(\tau)\leq -a_0b_0/2$ for every $\tau\in[-\tau_0,\tau_0]$. We have to distinguish three cases:
\begin{itemize}
 \item $F(\tau_1)\leq 0$. In that case, we take $\tau'=\tau_1$ and $\tau=\tau_2$ and we find
 $$r(\tau_2)-r(\tau_1)\leq-\frac{a_0b_0}{4}(\tau_2-\tau_1)^2+\mathcal{O}((\tau_2-\tau_1)^3).$$
 From this, we can deduce that $|\tau_2-\tau_1|\leq \frac{32r^{1/2}}{a_0b_0}.$
 \item $F(\tau_1)\geq 0$ and $F(\tau_2)\geq 0$. In that case, we take $\tau'=\tau_2$ and $\tau=\tau_1$ and we find
 $$r(\tau_1)-r(\tau_2)\leq-\frac{a_0b_0}{4}(\tau_2-\tau_1)^2+\mathcal{O}((\tau_2-\tau_1)^3).$$
 Again, we deduce an upper bound of order $\mathcal{O}(r^{\frac{1}{2}}).$
 \item $F(\tau_1)>0$ and $F(\tau_2)<0$. In that case, one can find some $\tau_0\in[\tau_1,\tau_2]$ such that $F(\tau_0)=0$. Then, we apply the above inequality twice to get 
 $$r(\tau_2)-r(\tau_0)=\frac{F'(\tau_0)}{2}(\tau_2-\tau_0)^2+\mathcal{O}((\tau_2-\tau_0)^3)\ \text{and}\ r(\tau_1)-r(\tau_0)=\frac{F'(\tau_0)}{2}(\tau_1-\tau_0)^2+\mathcal{O}((\tau_1-\tau_0)^3).$$
 Combining the two equalities, we find that $|\tau_2-\tau_1|=\mathcal{O}(r^{\frac{1}{2}}).$
\end{itemize}
Gathering these bounds, we find that, for every $x\in B_{\epsilon_0}(x_0)$ and for every $r\leq r_1$,
\begin{equation}\label{e:bound-tangent-case}\forall\gamma\in\mathcal{U}_{\epsilon_0},\ 0\leq\frac{1}{2\tau_0}\int_{-\tau_0}^{\tau_0}\mathcal{R}(\mathbf{1}_{B_{2r}(x)})\circ\varphi_{\langle V\rangle}^{\tau}(\gamma)d\tau\leq Cr^{\frac{3}{2}}.\end{equation}

\subsubsection{The conclusion} 

By compactness, one can find $\gamma_1,\ldots,\gamma_N$ in $\Gamma_{x_0}$ and $\epsilon_1,\ldots\epsilon_N>0$ such that $\cup_{j=1}^N\mathcal{U}_{\epsilon_j}(\gamma_j)$ covers $\Gamma_{x_0}$. We take $\epsilon_0:=\min\{\epsilon_j:1\leq j\leq N\}$ so that $\mathcal{A}_{2\epsilon_0}\subset\cup_{j=1}^N\mathcal{U}_{\epsilon_j}(\gamma_j).$ In particular, given any $x\in B_{\epsilon_0}(x_0)$ and any $r<r_1$ (with $r_1$ chosen small enough to handle each neighborhood $\mathcal{U}_{\epsilon_j}(\gamma_j)$), the support of the map 
$$\gamma\mapsto \frac{1}{2\tau_0}\int_{-\tau_0}^{\tau_0}\mathcal{R}(\mathbf{1}_{B_{2r}(x)})\circ\varphi_{\langle V\rangle}^{\tau}(\gamma)d\tau$$
is contained in $\cup_{j=1}^N\mathcal{U}_{\epsilon_j}(\gamma_j)$. Thus, applying~\eqref{e:bound-transverse-case} and~\eqref{e:bound-tangent-case} to~\eqref{e:final-bound-ball}, we obtain, for any normalized solution $u_h$ to~\eqref{e:semiclassical-schrodinger}, 
$$\int_{\IS^2}\chi_{x,r}(y)|u_h(y)|^2d\upsilon_{g_0}(y)=\mathcal{O}(r^{\frac{3}{2}})+\mathcal{O}(h^{1-3\beta}),$$
where the constant can be made uniform for $x\in B_{\epsilon_0}(x_0)$ and $r\geq h^{\beta}$. Taking $\beta= \frac{2}{9}$ yields Proposition~\ref{p:main-prop}.

\begin{rema}\label{r:caustic} The analysis of the vector field performed here is related to the analysis in~\cite[\S~4]{MR16}. In that reference, we showed with Maci\`a that the semiclassical measures of $-\Delta_{g_0}+V$ can be decomposed as a convex combination of the Haar measures carried by the Lagrangian tori of the completely integrable system $(H_0,\mathcal{R}_{\text{cl}}(V))$. For $2$-dimensional tori, the projection of the Haar measure on $\IS^2$ is absolutely continuous~\cite[Th.~4.3]{MR16} with some eventual blow-up of the density at some points which are often called caustics~\cite[Lemma~4.6]{MR16}. This regularity of the projection is exactly the property we have been using here in a somewhat refined way to get our bounds $\mathcal{O}(r^{1+\alpha})$. The bound~\eqref{e:bound-transverse-case} ($\alpha=1$) corresponds to points of these $2$-dimensional Lagrangian tori where the projection is regular while~\eqref{e:bound-tangent-case} ($\alpha=1/2$) corresponds to these caustics. 
\end{rema}

%As above, for $i=1,2$ (but now at order $2$), we write
%$$\phi(\tau_{2i})-\phi(\tau_{2i-1})=(\tau_{2i}-\tau_{2i-1})\phi'(\tau_{2i-1})+\frac{1}{2}(\tau_{2i}-\tau_{2i-1})^2\phi''(\tau_{2i-1})+\mathcal{O}((\tau_{2i}-\tau_{2i-1})^3),$$
%where the constant in the remainder can be made uniform in terms of $r$, $\gamma$ and $\tau_j$. Again, one has
%$$|\phi(\tau_{2i})-\phi(\tau_{2i-1})|\leq 4r+\left|\text{arccos}\left(-\cos(\theta(\tau_{2i-1})-\theta_x)\sin(\phi_x)\right)-\text{arccos}\left(-\cos(\theta(\tau_{2i})-\theta_x)\sin(\phi_x)\right)\right|.$$

%Contrary to the transverse case, we rather make use of the equation~\eqref{e:HJ} for $\theta$ to deduce that, for every $-\tau_0\leq\tau_1\leq\tau_2\leq\tau_0$, 
%\begin{equation}\label{e:angle-interval-transverse}\theta(\tau_2)-\theta(\tau_1)\geq a_0(\tau_2-\tau_1).\end{equation}

%We will now verify that there exists some $r_0, r_1>0$ and $c_0>0$ such that, for every $\gamma\in U_{\gamma_0}$ and for every $0<r<r_1$
%\begin{equation}\label{e:travel-time}
%\forall x\in B_{r_0}(x_0),\quad \text{Leb}\left(\left\{\tau\in[-\tau_0/4,\tau_0/4]: \mathcal{R}(\mathbf{1}_{B_{2r}(x)})\circ\varphi_{\langle V\rangle}^{\tau}(\gamma)\neq 0\right\}\right)\leq c_0 r.
%\end{equation}
%To see this, we can proceed by contradiction and suppose that, for every $n\geq 1$, one can find $\gamma_n \in U_{\gamma_0}$, $0<r_n<\frac{1}{n}$ and $x_n\in B_{1/n}(x_0)$ so that
%$$\text{Leb}\left(\left\{\tau\in[-\tau_0/4,\tau_0/4]: \mathcal{R}(\mathbf{1}_{B_{2r_n}(x_n)})\circ\varphi_{\langle V\rangle}^{\tau}(\gamma_n)\neq 0\right\}\right)\geq nr_n.$$
 
 \section{Final comments}

\subsection{Relaxing assumption~\eqref{e:hyp-trans}} 
 
Up to some extra work, assumption~\eqref{e:hyp-trans} could certainly be relaxed. For instance, one could require instead that the critical points are of finite order i.e. the derivative does not vanish at a certain order which may be larger than $2$. We would then end up with some upper bound of order $\mathcal{O}(r^{1+\alpha})$ for some $0<\alpha\leq 1/2$ related to the order of vanishing at the critical points of $\mathcal{R}(V)|_{\Gamma_{x_0}}.$ This would give slightly worst upper bound on the growth of $L^p$-norms but it would allow to take larger compact subsets $K$ in~\eqref{e:Lp-compact}.

 \subsection{Relaxing assumption~\eqref{e:hyp-crit}} 
A priori, it does not seem possible to remove assumption~\eqref{e:hyp-crit} from the hypothesis of Proposition~\ref{p:main-prop}. Indeed, if there exists $\gamma_0\in\Gamma_{x_0}$ such that $X_{\langle V\rangle}(\gamma_0)=0$, then the value of~\eqref{e:av-radon-potential} at $\gamma_0$ will be equal to $4r$ and it will prevent us from drawing the same conclusion using our argument.

\subsection{Sharpness of the exponents} Even if we tried to optimize our arguments, it is not clear if the bounds we obtain on $L^p$-norms are sharp or not. The argument works as well for elements in $L^2(\IS^2)$ which verify~\eqref{e:schrodinger} modulo some small remainder (say $\ml{O}(\lambda^{-2})$) and it would be interesting (but probably subtle) to construct quasimodes saturating these local $L^p$-estimates.

\subsection{The range $p>6$}
In this range, it is plausible that the methods from~\cite{GaTo18, Ga19, CaGa19, CaGa20} allow to handle these critical geodesics. Indeed, suppose that there exist a point $x_0\in\mathbb{S}^2$ and a sequence $(\psi_{\lambda_k})_{k\geq 1}$ of normalized solutions to~\eqref{e:schrodinger} verifying $\lambda_k\rightarrow+\infty$ and
\begin{equation}\label{e:blowup}\lim_{k\rightarrow+\infty}|\psi_{\lambda_k}(x_0)|\lambda_k^{-\frac{1}{2}}\neq0.\end{equation}
Up to extracting a subsequence, we can suppose that $(\psi_{\lambda_k})_{k\geq 1}$ has a single semiclassical measure $\mu$~\cite[Ch.~5]{Zw12}. Recall that it is a probability measure carried by $S^*\IS^2$ which is invariant by the geodesic flow $\varphi_0^t$. In particular, it induces a measure $\tilde{\mu}$ on $G(\IS^2)$. Then, we can consider $\tilde{\mu}_{x_0}=\tilde{\mu}|_{\Gamma_{x_0}}$. This measure can be decomposed into three parts: the absolutely continuous component, the singular continuous one and the pure point one. According to the results of Galkowski and Toth in~\cite{GaTo18}, property~\eqref{e:blowup} implies that the absolutely continuous part is not identically $0$. Combined with~\cite[Prop.~2.3]{MR16}, this implies that $\mathcal{R}(V)|_{\Gamma_{x_0}}$ has infinitely many critical points. In other words, if $\mathcal{R}(V)|_{\Gamma_{x_0}}$ has finitely many critical points, then, for any sequence $(\psi_{\lambda_k})_{k\geq 1}$ of normalized solutions to~\eqref{e:schrodinger}, one has 
$$|\psi_{\lambda_k}(x_0)|=o\left(\lambda_k^{\frac{1}{2}}\right),$$
which improves the remainder from the local Weyl law at $x_0$ without imposing~\eqref{e:hyp-crit}. Compared with Theorem~\ref{t:maintheo}, this is of course not quantitative. If one is able to combine the quantitative arguments of Canzani and Galkowski~\cite{CaGa19, CaGa20} with the extra invariance by the flow of $X_{\langle V\rangle}$~\cite{MR16}, then this may give rise to improvements on Sogge's upper bounds~\eqref{e:Lp-sogge} in the range $p>6$ under weaker geometric assumptions than the ones appearing in Theorem~\ref{t:maintheo}. Recall from the introduction that, thanks to the conjugation formula~\eqref{e:guillemin-weinstein}, eigenfunctions of $-\Delta_{g_0}+V$ which are the image under $\mathcal{U}$ of joint eigenfunctions for $(-\Delta_{g_0},V^\sharp)$ enjoy improved $L^p$ estimates near $x_0$ (for $p>6$) under appropriate assumptions on the critical points of $\mathcal{R}(V)|_{\Gamma_{x_0}}$~\cite{GaTo20, Ta19b}. In particular, if the spectrum of $-\Delta_{g_0}+V$ is simple~\cite[Th.~7]{Uh76}, then all eigenfunctions of $-\Delta_{g_0}+V$ will be the image of joint eigenfunctions.

 \subsection{The case of odd potentials}
 
 In~\cite{MR19}, it was shown that one can uncover extra-invariance properties of semiclassical measures even if $\mathcal{R}(V)$ identically vanishes (meaning that $V$ is an odd function, e.g. $V(x_1,x_2,x_3)=x_3$). In principle, the above arguments could be adapted following the lines of this reference, up to some extra technical work. In that case, the role of $\mathcal{R}(V)$ would be played by the function
$$\mathcal{R}^{(2)}(V)=\mathcal{R}(V^2)-\frac{1}{2\pi}\int_0^{2\pi}\int_0^t\{V\circ\varphi_0^t,V\circ\varphi_0^s\}dsdt.$$
See also~\cite{Gu78, Ur85} for earlier related results on spectral asymptotics of Schr\"odinger operators.

\subsection{Semiclassical operators}

In Remarks~\ref{r:Lp-semiclassical} and~\ref{r:KN-Lp}, we observed that our bounds on $L^p$ norms are valid more generally for solutions to
$$-h^2\Delta_{g_0}u_h+\varepsilon_h Vu_h=u_h,\quad\|u_h\|_{L^2(\IS^2)}=1.$$
Even if it was maybe not optimal, for $p>6$, we needed to impose $\varepsilon_h\leq h^{1+\epsilon}$ for some positive $\epsilon$ while for $4\leq p<6$, we only required $\varepsilon_h\leq h$. Thanks to Remarks~\ref{r:semiclassical-pert} and~\ref{r:semiclassical-pert-2}, this yields the following bounds on $L^p$ norms. For $p=\infty$, one has
$$\|u_h\|_{L^{\infty}(B_{r_0}(x_0))}\leq C_{\infty,x_0} h^{-\frac{1}{2}}\left(h^{\frac{1}{18}}+h^{\frac{\epsilon}{4}}\right),$$
which yields a polynomial improvement over the usual bound. In the range $4<p<6$, we get similarly, for any $r\geq h^{\frac{2}{9}}$,
$$\|u_h\|_{L^{p}(B_{r_0}(x_0))}\leq C_{p,x_0} h^{-\sigma_0(p)}\left(h^{\frac{1}{9}}+(rh)^{-1}\varepsilon_h\right)^{\frac{1}{2}\left(\frac{6}{p}-1\right)},$$
while for $p=4$, we end up with
$$\|u_h\|_{L^{4}(B_{r_0}(x_0))}\leq C_{4,x_0} |\log h|h^{-\frac{1}{8}}\left(h^{\frac{1}{9}}+(rh)^{-1}\varepsilon_h\right)^{\frac{1}{4}},$$
In these last two cases, it yields improvements over Sogge's upper bound as soon as $h^{-1}\varepsilon_h\rightarrow 0$. Note that in every cases, $\varepsilon_h$ may go to $0$ very fast. For instance, one may have $\varepsilon_h\ll h^2$. 

\subsection{The case of Zoll surfaces}
Following the lines of~\cite{MR16}, we could adapt the results to Laplace eigenfunctions,
$$-\Delta_g\psi_\lambda=\lambda^2\psi_\lambda,$$
where $g$ is a $C_{2\pi}$ (or Zoll) metric on $\IS^2$, i.e. all of whose geodesics are closed, simple and of length $2\pi$. See~\cite{Bes78} for a detailed review on this geometric assumption. In that case, it is known~\cite{CdV79} that
$$\sqrt{-\Delta_g}=A+\frac{\alpha}{4}+Q,$$
where $Q$ is a pseudodifferential operator of order $-1$, $\alpha$ is the Maslov index of the closed trajectories and $\text{Sp}(A)\subset\IZ_+$. Combining the above proof with the arguments from~\cite[\S 3.1]{MR16}, we will end up with the same quantities as in~\eqref{e:final-bound-ball} except that $\mathcal{R}(V)$ will be replaced by some function $q_0(x,\xi)$ (related to the principal symbol of $Q$). An exact expression for $q_0$ was given by Zelditch in~\cite{Ze96, Ze97} and it involves curvature terms of the metric. Under the geometric assumptions of Theorem~\ref{t:maintheo} on the point $x_0$ but with $q_0$ replacing $\mathcal{R}(V)$, we could obtain improved $L^p$-bounds near $x_0$. Yet, the expression of $q_0$ being a little bit involved, this condition is harder to verify.

\subsection{The higher dimensional case}

For the sake of simplicity, we restricted ourselves to the $2$-dimensional case but the extra invariance property by the flow of $X_{\langle V\rangle}$ remains true in higher dimensions $n\geq 3$~\cite[Prop.~2.3]{MR16}. Thus, modulo some extra work and some appropriate assumptions on $X_{\langle V\rangle}|_{\Gamma_{x_0}}$, one should be able to obtain localized $L^2$-estimates as in Proposition~\ref{p:main-prop} but maybe for smaller values of $\alpha$. Then, in the range $p_c=\frac{2(n+1)}{n-1}<p\leq+\infty$, this can be transferred into $L^p$ bounds using that~\eqref{e:localized-Lp} remains true for $p=\infty$ in dimension $n\geq 3$~\cite[Eq.(3.3)]{So16}. Similarly, for $p<p_c$, the Kakeya-Nikodym bounds of Section~\ref{s:kakeya} remains true up to $p>\frac{2(n+2)}{n}$ and they can again be roughly bounded by the $L^2$-localized norms appearing in Proposition~\ref{p:main-prop}. Yet, we are not aware of an analogue of Guillemin's Theorem~\cite{Gu76} showing that $\mathcal{R}$ is an isomorphism when restricted to the appropriate spaces of smooth functions on $\IS^n$ and $G(\IS^n)$ and hence making the condition on $x_0$ easy to verify.


\begin{thebibliography}{99}
\bibitem{Be77} P.B\'erard, \emph{On the wave equation on a compact Riemannian manifold without conjugate points}, Math. Z. $\mathbf{155}$ (1977), 249--276
\bibitem{Bes78} A. Besse, \emph{Manifolds All of Whose Geodesics Are Closed}, Ergeb. Math. 93, Springer-Verlag, New York (1978)
\bibitem{BlSo15} M. Blair and C.~Sogge, \emph{Refined and microlocal Kakeya-Nikodym bounds for eigenfunctions in two dimensions}, Analysis and PDE $\mathbf{8}$ (2015), 747--764
\bibitem{BlSo17} M. Blair and C.~Sogge, \emph{Refined and Microlocal Kakeya-Nikodym Bounds of Eigenfunctions in Higher Dimensions}, Comm. in Math. Phys. $\mathbf{356}$ (2017), 501--533
\bibitem{BlSo18} M. Blair and C.~Sogge, \emph{Concerning Toponogov's Theorem and logarithmic improvement of estimates of eigenfunctions}, Journal of Differential Geometry, $\mathbf{109}$ (2018), 189--221.
\bibitem{BlSo19} M. Blair and C.~Sogge, \emph{Logarithmic improvements in $L^p$ bounds for eigenfunctions at the critical exponent in the presence of nonpositive curvature}, Inv. math. $\mathbf{217}$ (2019), 703--748
\bibitem{Bon17} Y.~Bonthonneau, \emph{The $\theta$ Function and the Weyl Law on Manifolds Without Conjugate Points}, Doc. Math. $\mathbf{22}$ (2017), 1275--1283
\bibitem{Bo93} J. Bourgain, \emph{Eigenfunctions  bounds  for  the  Laplacian  on  the $n$-torus}, IMRN (1993), 61--66
\bibitem{Bo09} J.~Bourgain, \emph{Geodesic restrictions and $L^p$-estimates for eigenfunctions of Riemannian surfaces}, in Linear and complex analysis, Amer. Math. Soc. Transl. $\mathbf{226}$ (2009), Amer. Math. Soc.,Providence, RI, 27--35.
\bibitem{BoDe15}J. Bourgain and C. Demeter, \emph{The proof of the $l^2$ decoupling conjecture}, Ann. of Math. $\mathbf{182}$ (2015), 351--389
\bibitem{BrLM20} S.~Brooks and E.~Le Masson, \emph{$L^p$ norms of eigenfunctions on regular graphs and on the sphere}, IMRN (2020), 3201--3228
\bibitem{BuLe13} N. Burq and G. Lebeau, \emph{Injections de Sobolev probabilistes et applications}, Ann. Sci. \'Ec. Norm. Supe\'er. $\mathbf{46}$ (2013), 917--962
\bibitem{CaGa19} Y.~Canzani and J.~Galkowski, \emph{Eigenfunction concentration via geodesic beams}, preprint arXiv1903.08461, to appear in J. Reine Angew. Math. (2019)
\bibitem{CaGa20} Y.~Canzani and J.~Galkowski, \emph{Growth of high $L^p$ norms for eigenfunctions: an application of geodesic beams}, preprint arXiv2003.02525 (2020)
\bibitem{CM11} T. Colding and W. P. Minicozzi II, \emph{Lower bounds for nodal sets of eigenfunctions}, Comm.Math. Phys. $\mathbf{306}$ (2011), 777--784
\bibitem{CdV79} Y.~Colin de Verdi\`ere, \emph{Sur le spectre des op\'erateurs elliptiques \`a  bicaract\'eristiques toutes p\'eriodiques}, Comment. Math. Helv. $\mathbf{54}$ (1979), 508--522
\bibitem{Co71} R. Cooke, \emph{A Cantor-Lebesgue Theorem in two dimensions}, Proc. AMS $\mathbf{30}$ (1971), 547--550
\bibitem{DJN19} S.~Dyatlov, L.~Jin and S.~Nonnenmacher, \emph{Control of eigenfunctions on surfaces of variable curvature}, preprint arXiv:1906.08923, to appear in JAMS (2019)
\bibitem{Ga19} J.~Galkowski, \emph{Defect measures of eigenfunctions with maximal $L^\infty$ growth}, Annales de L'institut Fourier $\mathbf{69}$ (2019), 1757--1798
\bibitem{GaTo18} J.~Galkowski and J.~Toth, \emph{Eigenfunction scarring and improvements in $L^\infty$ growth}, Anal. PDE $\mathbf{11}$ (2018), 801--812
\bibitem{GaTo20} J.~Galkowski and J.~Toth, \emph{Pointwise bounds for joint eigenfunctions of quantum completely integrable systems}, Comm. Math. Phys. $\mathbf{375}$ (2020), 915--947
\bibitem{Gu76} V. Guillemin, \emph{The Radon transform on Zoll surfaces}, Adv. Math. $\mathbf{22}$ (1976), 85--119
\bibitem{Gu78} V.~Guillemin, \emph{Some spectral results for the Laplace operator with potential on the $n$-sphere}, Adv. in Math. $\mathbf{27}$ (1978), 273--286
\bibitem{Gu81} V.~Guillemin, \emph{Band asymptotics in two dimensions}, Adv. in Math. $\mathbf{42}$ (1981), 248--282
\bibitem{HaTa15} A.~Hassell and M.~Tacy, \emph{Improvement of eigenfunction estimates on manifolds of nonpositive curvature}, Forum Mathematicum $\mathbf{27}$ (2015), 1435--1451
\bibitem{HR16} H.~Hezari and G.~Rivi\`ere, \emph{$L^p$ norms, nodal sets and quantum ergodicity}, Adv.Math. $\mathbf{290}$ (2016), 938--966
\bibitem{Ho68} L.~H\"ormander, \emph{The spectral function of an elliptic operator}, Acta  Math. $\mathbf{121}$ (1968), 193--218
\bibitem{IwSar95} H. Iwaniec and P. Sarnak, \emph{$L^{\infty}$-norms of eigenfunctions of arithmetic surfaces}, Ann. Math. $\mathbf{141}$ (1995), 301--320
\bibitem{KTZ} H.~Koch, D.~Tataru and M.~Zworski, \emph{Semiclassical $L^p$ estimates}, Ann. H. Poincar\'e $\mathbf{8}$ (2007), 885--916
\bibitem{MR16} F.~Maci\`a and G.~Rivi\`ere, \emph{Concentration and Non-Concentration for the Schr\"odinger Evolution on Zoll Manifolds}, Comm. Math. Phys. $\textbf{345}$ (2016), 1019--1054
\bibitem{MR19} F.~Maci\`a and G.~Rivi\`ere, \emph{Observability and quantum limits for the Schrödinger equation on the sphere}, in \emph{Probabilistic Methods in Geometry, Topology and Spectral Theory}, Contemporary Mathematics $\mathbf{739}$ (2019)
\bibitem{Ma16} S.~Marshall, \emph{$L^p$ norms of higher rank eigenfunctions and bounds for spherical functions}, J. Eur. Math. Soc. $\mathbf{18}$ (2016), 1437--1493
\bibitem{Sa88} Y. Safarov, \emph{Asymptotics of a spectral function of a positive elliptic operator without a non trapping condition}, Funktsional. Anal. i Prilozhen 22 (1988), 53--65, translated in Funct. Anal. Appl. 22 (1988), 213--223
\bibitem{Sar04} P.~Sarnak, Letter to Morawetz, Available at http://publications.ias.edu/sarnak (2004)
\bibitem{Sh74} A.I.~Shnirelman, \emph{Ergodic properties of eigenfunctions}, Uspehi Mat. Nauk $\mathbf{180}$ (1974), 181--182
\bibitem{Sh74b} A.I.~Shnirelman, \emph{Statistical properties of eigenfunctions}, In Proceedings  of  the  All-USSR  School  in Differential Equations with Infinite Number of Independent Variables and in Dynamical Systems with Infinitely Many Degrees of Freedom, May 1973. Armenian Academy of Sciences, Erevan, 1974. Translation avalaible at http://math.mit.edu/$\sim$dyatlov/files/2019/shnirelman.pdf.
\bibitem{So88} C.~Sogge, \emph{Concerning the $L^p$ norm of spectral clusters for second-order elliptic operators on compact manifolds}, J. Funct. Anal. $\mathbf{77}$ (1988) 123--138
\bibitem{So11} C.~Sogge, \emph{Kakeya-Nikodym  averages  and $L^p$-norms of eigenfunctions}, Tohoku  Math. J. $\mathbf{63}$ (2011), 519--538
\bibitem{So15} C.~Sogge, \emph{Problems related to the concentration of eigenfunctions}, Journ\'es \'equations aux d\'eriv\'ees partielles (2015), avalaible at http://www.numdam.org/
\bibitem{So16} C.~Sogge, \emph{Localized $L^p$-estimates of eigenfunctions: a note on an article of Hezari and Rivi\`ere}, Adv. Math. $\mathbf{289}$ (2016), 384--396
\bibitem{So17} C.~Sogge, \emph{Improved  critical  eigenfunction  estimates  on  manifolds  of  nonpositive curvature}, Math. Res. Lett. $\mathbf{24}$ (2017), 549--570
\bibitem{SoToZe11} C. Sogge, J.A.~Toth, and S.~Zelditch, \emph{About the blowup of quasimodes on Riemannian manifolds}, J. Geom. Anal. $\mathbf{21}$ (2011),  150--173.
\bibitem{SoZe02} C.~Sogge and S.~Zelditch, \emph{Riemannian manifolds with maximal eigenfunction growth}, Duke Math. J. $\mathbf{114}$ (2002), 387--437
\bibitem{SoZe16} C.~Sogge and S.~Zelditch, \emph{Focal points and sup-norms of eigenfunctions}, Rev. Mat. Iberoam. $\mathbf{32}$ (2016), 971--994
\bibitem{Ta19} M.~Tacy, \emph{$L^p$ estimates for joint quasimodes of semiclassical pseudodifferential operators}, Israel Journal of Mathematics $\mathbf{232}$ (2019), 401--425
\bibitem{Ta19b} M.~Tacy, \emph{$L^p$ estimates for joint quasimodes of semiclassical pseudodifferential operators whose characteristic sets have kth order contact}, preprint arXiv:1909.12559 (2019)
\bibitem{To96} J.~Toth, \emph{Eigenfunction localization in the quantized rigid body}, J. Differential Geom. $\mathbf{43}$ (1996), 844--858
\bibitem{ToZe02}  J.~Toth and S.~Zelditch, \emph{Riemannian manifolds with uniformly bounded eigenfunctions}, Duke Math. J. $\mathbf{111}$ (2002), 97--132
\bibitem{ToZe03} J.~Toth and S.~Zelditch, \emph{$L^p$ norms of eigenfunctions in the completely integrable case}, Ann. Henri Poincar\'e $\mathbf{4}$ (2003), 343--368
\bibitem{Uh76} K. Uhlenbeck, \emph{Generic properties of eigenfunctions}, American J. Math. $\mathbf{98}$ (1976), 1059--1078
\bibitem{Ur85} A.~Uribe, \emph{Band invariants and closed trajectories on $\IS^n$}, Adv. in Math. $\mathbf{58}$ (1985), 285--299
\bibitem{VdK97} J.M.~Van der Kam, \emph{$L^{\infty}$-norms and quantum ergodicity on the sphere}, IMRN (1997), 329--347
\bibitem{Wa11} W.M.~Wang, \emph{Eigenfunction Localization for the 2D Periodic Schrodinger Operator}, IMRN (2011), 1804--1838
\bibitem{W77} A.~Weinstein, \emph{Asymptotics of eigenvalue clusters for the Laplacian plus a potential}, Duke Math. Jour. $\mathbf{44}$ (1977), 883--892
\bibitem{Ze96} S. Zelditch, \emph{Maximally degenerate Laplacians}, Ann. Inst. Fourier $\mathbf{46}$ (1996), 547--587
\bibitem{Ze97} S.~Zelditch, \emph{Fine structure of Zoll spectra}, J. Funct. Anal. $\mathbf{143}$ (1997), 415--460
\bibitem{Ze08} S.~Zelditch, \emph{Local and global analysis of eigenfunctions}, Advanced Lectures in Mathematics $\mathbf{7}$ (2008), 545--658
\bibitem{Ze17} S.~Zelditch, \emph{Eigenfunctions of the Laplacian on a Riemannian manifold}, CBMS Regional Conference Series in Mathematics, 125. Published for the Conference Board of the Mathematical Sciences, Washington, DC; by the American Mathematical Society, Providence, RI, 2017
\bibitem{Zh20} Y.~Zhang, \emph{On Fourier restriction type problems on compact Lie groups}, preprint arXiv:2005.11451 (2020)
\bibitem{Zh21} Y.~Zhang, \emph{Schr\"odinger equations on compact globally symmetric spaces}, preprint arXiv:2005.00429, to appear in J. Geom. Anal. (2021)
\bibitem{Zy74} A. Zygmund, \emph{On Fourier coefficients and transforms of functions of two variables}, Studia Mathematica $\mathbf{50}$ (1974), 189--201 
\bibitem{Zw12} M.~Zworski \emph{Semiclassical analysis}, Graduate Studies in Mathematics $\mathbf{138}$, AMS (2012)
\end{thebibliography}
\end{document}